\documentclass[12pt,reqno]{amsart}
\usepackage{amsmath,amsfonts,amssymb,amscd,amsthm,amsbsy,epsf}
\usepackage[dvips]{graphicx}
\usepackage[usenames]{xcolor}
\usepackage{hyperref}
\usepackage{comment}

\numberwithin{equation}{section}
\newtheorem{theorem}{Theorem}
\newtheorem{meta-thm}[theorem]{Meta-Theorem}

\newtheorem{proposition}[theorem]{Proposition}

\DeclareMathAlphabet{\mathcalligra}{T1}{calligra}{m}{n}

\newcommand\beq[1]{ \begin{equation}\label{#1} }
\newcommand{\eeq}{ \end{equation} }

\newcommand\beqa[1]{ \begin{eqnarray} \label{#1}}

\newcommand{\eeqa}{ \end{eqnarray} }

\newcommand{\beqano}{ \begin{eqnarray*} }

\newcommand{\eeqano}{ \end{eqnarray*} }

\newcommand\equ[1]{{\rm (\ref{#1})}}

\newcommand{\nn}{\nonumber \\}

\def\C{{\mathcal C}}

\def\P{{\mathcal P}}
\def\E{{\mathcal E}}

\def\H{{\mathcal H}}
\def\K{{\mathcal K}}
\def\P{{\mathcal P}}
\def\R{{\mathcal R}}

\def\Z{{\mathcal Z}}

\def\integer{{\mathbb Z}}
\def\nat{{\mathbb N}}

\def\real{{\mathbb R}}

\begin{document}
\title[The dynamics around the collinear points]
{The dynamics around the collinear points of the elliptic 
three-body problem: a normal form approach}

\author[A. Celletti]{Alessandra Celletti}
\address{
Department of Mathematics, University of Roma Tor Vergata, Via
della Ricerca Scientifica 1, 00133 Roma (Italy)}
\email{celletti@mat.uniroma2.it}

\author[C. Lhotka]{Christoph Lhotka}
\address{
Department of Mathematics, University of Roma Tor Vergata, Via
della Ricerca Scientifica 1, 00133 Roma (Italy)}
\email{lhotka@mat.uniroma2.it}

\author[G. Pucacco]{Giuseppe Pucacco}
\address{
Department of Physics, University of Roma Tor Vergata, Via della
Ricerca Scientifica 1, 00133 Roma (Italy)}
\email{pucacco@roma2.infn.it}

\baselineskip=18pt              


\maketitle

\begin{abstract}
We study the dynamics of the collinear points in the planar, restricted  three-body problem, assuming that the primaries move on an elliptic orbit around a common barycenter. The equations of motion can be conveniently written in a rotating-pulsating barycentric frame, taking the true anomaly as independent variable. We consider the Hamiltonian modeling this problem in the extended phase space and we implement a normal form to make a center manifold reduction. The normal form provides an approximate solution for the Cartesian coordinates, which allows us to construct several kinds of orbits, most notably planar and vertical Lyapunov orbits, and halo orbits. We compare the analytical results with a numerical simulation, which requires special care in the selection of the initial conditions. 
\end{abstract}



\section{Introduction}
Despite their instability, the collinear equilibrium points of the planar, circular, restricted three-body problem (hereafter, PCR3BP) play a key role in Celestial Mechanics and Astrodynamics. In the PCR3BP, a small particle moves in the gravitational field of two primaries, which orbit around each other on a circular orbit; moreover, all bodies move in the same plane. The collinear points, which are equilibria in the synodic reference frame, have been used to design space trajectories (\cite{Conley}, \cite{kolomaro}) and they are considered  privileged positions where to place spacecraft or instruments. There are many works devoted to analytical and numerical solutions around the collinear points in the circular problem (see, e.g., \cite{JM}, \cite{Richardson}, \cite{GM}, \cite{CPS}, \cite{CCP}). 
In this work, we extend the study of the dynamics around the collinear points assuming that the primaries move on an elliptic orbit, thus providing a more \sl realistic \rm description. 

An analytical approach based on the construction of a resonant normal form has 
been applied in \cite{CCP}, \cite{CPS} in the case of the spatial \sl circular \rm
problem. This approach provides an integrable 3-DOF Hamiltonian
which captures both the hyperbolic dynamics normal to the center
manifold and the elliptic dynamics on it. This analytical approach allows us also to compute
approximate solutions for periodic and quasi-periodic orbits,
including those associated to 1:1 resonant bifurcations. We remark that the method
can be easily generalised to include perturbing effects like the
Solar radiation pressure and the oblateness of the primaries
(\cite{AMPA}, \cite{LXD}). 

In this work, we extend the construction of the analytical solution to the elliptic case by computing a suitable normal form, which describes Lyapunov planar orbits, Lyapunov vertical orbits, halo orbits; these orbits are  characterized by prescriptions on the normal modes or by resonant conditions.  As it is well known (\cite{Broucke}, \cite{Szebehely}, \cite{Alebook}), in the planar, elliptic, restricted three-body problem, it is convenient to introduce a rotating-pulsating frame, which rotates with the true anomaly of the elliptic orbit, which is taken as independent variable; moreover, the coordinates are rescaled with the orbital radius. The Hamiltonian describing this problem depends on the true anomaly and therefore it is useful to consider it in the extended phase space by adding a dummy action conjugated to the true anomaly, thus obtaining a 4DOF Hamiltonian. 
In practice, however, the normal-form dynamics retains all the properties of the circular
case and it is effectively an integrable 3-DOF problem. The time
dependence is encoded in the generating function of the
normalising transformation and becomes apparent when the solutions are
mapped back to the original variables.

We provide the explicit expressions of the normal form Hamiltonian
for the collinear points $L_1$, $L_2$ and we give examples of the main families of periodic orbits
for several cases of interest in astrodynamical applications. We remark that, in
the work \cite{Holi}, the same problem
has been treated with the Poincar\'e-Lindstedt method. The
normal-form construction and the solutions thereof reproduces 
the results in \cite{Holi} with additional dynamical clues. A first order
resonant normal form already allows us to analytically compute the
bifurcation threshold of the halo family which, for $L_1$ and $L_2$, 
turns out to be slightly higher than in the corresponding
(namely, same mass ratio $\mu$) circular problem. 

In the specific case of the Earth-Moon system some 
Lyapounov and halo orbits are given together with some 
examples of isolated periodic orbits with higher commensurability 
(\cite{PeXu}, \cite{FeLa}). Explicit
expressions of the generating functions of the Lie transform are
given in order to compute invariant manifolds (periodic orbits and
invariant tori) at second order in the expansion.

The results are validated by a numerical integration of the equations of motion in which the initial conditions are given by the normal form procedure and carefully refined in order to find the given orbit, either Lyapunov or halo orbit. Finally, we mention that to get reliable trajectories, it is essential to use the symmetry of the periodic orbit which leads to some constraints in the initial conditions. 

\vskip.1in 

The paper is organized as follows. In Section~\ref{sec:CM} we give the setup in the elliptic case and we proceed with the center manifold reduction; in
Section~\ref{sec:NF1} we describe the resonant normal form; in Section~\ref{sec:halo} the normal form is used to
compute the bifurcation thresholds of the halo orbits; 
Section~\ref{sec:EMC} provides examples of different orbits 
for the $L_1$ equilibrium point in the Earth-Moon case; 
Section~\ref{sec:numerical} presents numerical simulations with a dedicated procedure to refine the initial conditions.

\section{Collinear points and center manifold reduction}\label{sec:CM}

In this Section, we present the Hamiltonian of the elliptic, restricted three-body problem in the rotating-pulsating frame (Section~\ref{sec:rotpuls}), we discuss the location of the collinear points (Section~\ref{sec:copo}), we introduce the Hamiltonian in the extended phase space (Section~\ref{sec:null}), we compute the quadratic part of the Hamiltonian in suitable coordinates (Section~\ref{sec:quadratic}) and we briefly recall the main results of the center manifold reduction (Section~\ref{normf}).

\subsection{The rotating-pulsating Hamiltonian}\label{sec:rotpuls}

We consider a satellite $S$ orbiting in the gravitational field of
two primaries, $P_1$ and $P_2$ with masses, respectively, $m_1$ and $m_2$ (with $m_2<m_1$).
We assume that the primaries orbit around their
common barycenter on a Keplerian ellipse with semimajor axis $a$
and eccentricity e. Denoting by $h$ the angular momentum constant,
related to the orbital elements by
$$
h^2=a(1-{\rm e}^2)\ ,
$$
the orbital radius  $r$ is given by
\beq{mp}
r(f)=\frac{h^2}{1+{\rm e} \cos f} \ ,
\eeq
where $f$ is the true anomaly.
 
We consider a \sl rotating-pulsating \rm barycentric frame, say $(O,X,Y,Z)$,
which is obtained as in \cite{Szebehely} (see also \cite{Alebook})
starting from a set of coordinates in a fixed frame with origin
coinciding with the barycenter $O$ of the primaries, transforming
the coordinates to a frame rotating with the true anomaly of the elliptic orbit
(rotating frame), then rescaling the coordinates by the radius $r$ (pulsating frame)
and taking the true anomaly, instead of time, as independent variable.
The latter choice implies that the time derivative is transformed into
$$
{d\over {dt}}={{df}\over {dt}}\ {d\over {df}}={h\over
{r(f)^2}}{d\over {df}}={{(1+{\rm e} \cos f)^2}\over {h^3}}{d\over
{df}}\ ;
$$
in the following, we will denote by a prime the derivative with
respect to the true anomaly.

In the non-uniformly rotating-pulsating frame $(O,X,Y,Z)$, the equations of
motion of the satellite are given by (see \cite{Szebehely},
Section 10.3.2.3):
\beqa{eq1}
X''-2Y'&=&{{\partial\Omega}\over {\partial X}}\nonumber\\
Y''+2X'&=&{{\partial\Omega}\over {\partial Y}}\nonumber\\
Z''+Z&=&{{\partial\Omega}\over {\partial Z}}\ ,
\eeqa
where the \sl potential \rm has the form
\beq{omegarp}
\Omega(X,Y,Z,f)={1\over {1+{\rm e}\cos f}}\ \left[{1\over
2}(X^2+Y^2+Z^2)+{{1-\mu}\over {r_1}}+{{\mu}\over {r_2}}+{1\over
2}\mu(1-\mu)\right]
\eeq
with 
$$\mu={m_2\over {m_1+m_2}}$$ 
and where the distances $r_1$, $r_2$ from the primaries, which are located at $(\mu,0,0)$, $(-1+\mu,0,0)$, are given by
$$
r_1=[(X-\mu)^2+Y^2+Z^2]^{1\over 2}\ ,\qquad  r_2=[(X+1-\mu)^2+Y^2+Z^2]^{1\over 2}\ .
$$
Denoting by $p_X$, $p_Y$, $p_Z$ the momenta conjugated to $X$,
$Y$, $Z$, the equations \equ{eq1} are associated to the 3D non-autonomous
Hamiltonian function
\beq{Hin}
H(p_X,p_Y,p_Z,X,Y,Z,f)={1\over
2}(p_X^2+p_Y^2+p_Z^2)+Yp_X-Xp_Y+{1\over
2}(X^2+Y^2+Z^2)-\Omega(X,Y,Z,f)\ .
\eeq
We consider such Hamiltonian on the collisionless manifold ${\mathcal
P}_c\subset\real^3\times\real^3$, defined as
$$
{\mathcal
P}_c\equiv\{(p_X,p_Y,p_Z),(X,Y,Z)\in\real^3\times\real^3:\
r_1(X,Y,Z)\not=0,\ r_2(X,Y,Z)\not=0\}
$$
and we endow ${\mathcal P}_c$ with the standard symplectic form
$$
\omega_s=dp_X\wedge dX+dp_Y\wedge dY+dp_Z\wedge dZ\ .
$$

\subsection{The collinear points}\label{sec:copo}

We denote by $L_1$, $L_2$, $L_3$ the collinear equilibrium points
inheriting the notation of the circular case. We underline that
the reference frames are oriented in such a way that $L_2$ lies at
the left of the smaller primary, while $L_3$ is at the right of the
larger primary.

Let $\gamma_j$, $j=1,2,3$, be the distances of the collinear
equilibrium points from the closest primary. 
The collinear equilibrium positions are the solutions of the equations
$$
{{\partial\Omega}\over {\partial X}}=0\ ,\qquad
{{\partial\Omega}\over {\partial Y}}=0\ ,\qquad
{{\partial\Omega}\over {\partial Z}}=0 
$$
with the condition $Y=Z=0$. It can be shown (see, e.g., \cite{Alebook}) that the
$\gamma_j$'s are the solutions of the following 5th order Euler's
equations:
\begin{equation}
\begin{split}
\gamma_1^5-(3-\mu)\gamma_1^4+(3-2\mu)\gamma_1^3-\mu\gamma_1^2+2\mu \gamma_1-\mu&=0  \qquad  \mbox{for}\hspace{1mm} L_1  \nonumber\\
\gamma_2^5+(3-\mu)\gamma_2^4+(3-2\mu)\gamma_2^3-\mu\gamma_2^2-2\mu \gamma_2-\mu&=0  \qquad   \mbox{for} \hspace{1mm} L_2 \nonumber\\
\gamma_3^5+(2+\mu)\gamma_3^4+(1+2\mu)\gamma_3^3-(1-\mu)\gamma_3^2-
2(1-\mu)\gamma_3-(1-\mu)&=0  \qquad \mbox{for} \hspace{1mm} L_3\ .
\end{split}
\end{equation}
As in the circular problem, for each $L_k, k=1,2,3$, we translate the origin so that it coincides with a collinear point: we introduce new coordinates $(x,y,z)$ through the following transformation
\beq{ctransf}
X= \mp \gamma_j x+\mu+a_j \ ,\qquad Y=\mp\gamma_j y \ ,\qquad Z=\gamma_j z \ .
\eeq
In \eqref{ctransf}, the upper signs hold
for $L_1$, $L_2$, while the lower signs hold for $L_3$ and
$a_1=-1+\gamma_1$, $a_2=-1-\gamma_2$, $a_3=\gamma_3$. In this way, the collinear `equilibria' in the elliptic problem, are actually pulsating with the rotating frame. 

The relation between the energy $E$ in rotating-pulsating synodic coordinates
and the energy $\widetilde{\mathcal E}$ in the local coordinates, in
virtue of the rescaling, is given by
\beqano
L_1: E &=&\widetilde{\mathcal E}\gamma_1^2-{1\over 2}(1-\gamma_1-\mu)^2-{\mu\over \gamma_1}-{{1-\mu}\over {1-\gamma_1}} \\
L_2: E &=&\widetilde{\mathcal E} \gamma_2^2-{1\over 2}(1+\gamma_2-\mu)^2-{\mu\over \gamma_2}-{{1-\mu}\over {1+\gamma_2}} \\
L_3: E &=&\widetilde{\mathcal E} \gamma_3^2-{1\over
2}(\gamma_3+\mu)^2-{{1-\mu}\over \gamma_3}-{{\mu}\over
{1+\gamma_3}}\ .
\eeqano
Introducing the coefficients
$c_n(\mu)$, $n\geq 2$, given by 
\beqano
c_n(\mu)&=&\frac{1}{\gamma_1^3}\left(\mu +(-1)^n \frac{(1-\mu)\gamma_1^{n+1}}{(1- \gamma_1)^{n+1}}\right) \qquad\hspace{2mm} \mbox{for} \hspace{2mm} L_1 \nonumber\\
c_n(\mu)&=&\frac{(-1)^n}{\gamma_2^3}\left(\mu+\frac{(1-\mu)\gamma_2^{n+1}}{(1+ \gamma_2)^{n+1}}\right) \qquad\hspace{6mm} \mbox{for} \hspace{2mm L_2} \nonumber\\
c_n(\mu)&=&\frac{(-1)^n}{\gamma_3^3}\left(1-\mu+\frac{\mu\gamma_3^{n+1}}{(1+
\gamma_3)^{n+1}}\right) \qquad \hspace{1mm}\mbox{for} \hspace{2mm}
L_3\
\eeqano
and using the Legendre polynomials $\P_n$, we expand the gravitational part of the potential \equ{omegarp} as
$$
-\frac{1-\mu}{r_1}-\frac{\mu}{r_2} = \sum_{n\geq 2}c_n(\mu)\rho^n
\P_n \left({x\over \rho}\right)\ , 
$$
where $\rho=(x^2+y^2+z^2)^{1\over 2}$. Introducing the functions
$$
T_n(x,y,z)=\rho^n\ \P_n\left({x\over \rho}\right)\ ,
$$
which are recursively defined through the formulae
$$
T_0=1\ , \qquad T_1=x\ ,\qquad
T_n=\frac{2n-1}{n}xT_{n-1}-\frac{n-1}{n}(x^2+y^2+z^2) T_{n-2} \qquad\forall n \ge 2 ,
$$
we obtain that the Hamiltonian \equ{Hin} can be written as
\beqa{ham1}
H(p_x,p_y,p_z,x,y,z,f)&=&\frac{1}{2}\left( p_x^2+p_y^2+p_z^2\right)+yp_x-xp_y +\nonumber \\
&& \frac12 \left(1-\frac{r(f)}{h^2}\right)(x^2+y^2+z^2)- \frac{r(f)}{h^2} \sum_{n\ge2}c_n(\mu) T_n (x,y,z) \ .\nonumber \\
\eeqa

\subsection{The Hamiltonian in the extended phase space}\label{sec:null}
The Hamiltonian \eqref{ham1} is 3D non-autonomous, due to the dependence on time of the true anomaly $f$. The standard
way to make the system autonomous is to \sl extend \rm the phase
space by considering $f$ as an angle in the Hamiltonian context
of action-angle variables and adding a new \sl dummy \rm action $F$, conjugated to $f$.
In the usual units in which $G(m_1+m_2)=1$ with $G$ the
gravitational constant and $a$=1, the time unit is such that the
revolution period of the primaries is $2 \pi$. Therefore the new
angle has unit frequency and, starting from \equ{ham1}, we can define the autonomous \sl null \rm Hamiltonian $\H$ as
\beq{ham2}
\H(p_x,p_y,p_z,F,x,y,z,f)= H + F \ ,
\eeq
where the dependence on $f$ appears through the radius $r(f)$.

Formally, we can introduce a new time $\tau$ which is equal to $f$; with a little abuse of notation,
we denote with a dot the derivative with respect to $\tau$ and write the equations of motion as
\beqano
\ddot x - 2 \dot y&=&{{\partial\Omega}\over {\partial x}},\nonumber\\
\ddot y +2 \dot x&=&{{\partial\Omega}\over {\partial y}},\nonumber\\
\ddot z +z&=&{{\partial\Omega}\over {\partial z}}\ ,
\eeqano
supplemented by
\beqano
\dot f&=&{{\partial\H}\over {\partial F}} = 1,\nonumber\\
\dot F&=&- {{\partial\H}\over {\partial f}} = {{\partial\Omega}\over {\partial f}}
= {{{\rm e} \sin f}\over {(1+{\rm e}\cos f})^2} W(x,y,z) \ ,
\eeqano
where $W$ denotes the effective potential
$$
W(x,y,z)={1\over
2}(x^2+y^2+z^2)+{{1-\mu}\over {r_1}}+{{\mu}\over {r_2}}+{1\over 2}\mu(1-\mu)\ .
$$
The Hamiltonian \equ{ham2} in the extended phase space 
$$ \H (p_x,p_y,p_z,F,x,y,z,f) = H + F \equiv 0
$$
is indeed conserved \rm and is always identically zero.
We remark that in the elliptic case there is no equivalent of the Jacobi constant that appears
in the circular case (see e.g. \cite{Szebehely}, Sect.10.3), therefore the conserved quantity \equ{ham2} for the extended system can be usefully exploited as a diagnostic in numerical computations.

\subsection{The quadratic part of the Hamiltonian}\label{sec:quadratic}
The general approach to the non-integrable system given by \eqref{ham1}-\eqref{ham2} is perturbative. 
Henceforth, we will write $\H$ as a series of the form
\beq{ham3}
\H(p_x,p_y,p_z,F,x,y,z,f)=\sum_{n\ge 0} \epsilon^n \H_n\ ,
\eeq
where $\epsilon$ is a \sl book-keeping \rm parameter used to
arrange terms in this and subsequent series. It is not necessarily related to a physical parameter,  
but has the role of providing a hierarchy of terms; in practice, 
after recursions usually implying multiplication of terms, 
a given series can be reordered and truncated and, in the end, $\epsilon$ can be set equal to one. 
A natural choice is that of considering of order $\epsilon^n$ the monomials in
$p_x,p_y,p_z,x,y,z$ of degree $n+2$. Moreover, 
since the dependence on $f$ is provided by the following series expansion
\beq{mpseries}
\frac{r(f)}{h^2}=\frac{1}{1+{\rm
e} \cos f} = \sum_{n=0}^{\infty}(-{\rm e})^n \cos^n f \ ,
\eeq
we can replace the 
expansion \eqref{mpseries} with the expression
\beq{mpsereps}
\frac{1}{1+{\rm e} \cos f} = \sum_{n=0}^{\infty} \epsilon^n (-{\rm
e})^n \cos^n f \ ,
\eeq
so that terms proportional to ${\rm e}^n$ are assumed to be of order $\epsilon^n$.

In order to put 
the terms of the starting Hamiltonian in a form most suitable for normalisation, we 
need to diagonalise the quadratic part. In view of the expansion \eqref{mpsereps}, already the
quadratic part  in
$p_x,p_y,p_z,x,y,z$ of the Hamiltonian
contains terms of order greater than zero in the book-keeping. 
We adopt the simplifying approach of keeping only truly zero-order terms in 
the quadratic part and therefore we perform 
the same diagonalising transformation of the
circular case on the quadratic Hamiltonian
\beq{ham2rq}
H_0^{(q)}(p_x,p_y,p_x,x,y,z)={1\over 2}(p_x^2+p_y^2+p_z^2)+yp_x-xp_y
-c_2x^2+{1\over 2}c_2 y^2+{1\over 2}c_2 z^2 \ .
\eeq 
Since the vertical component is already
diagonalised, one can limit to define a vector
$\xi\equiv(x,y,p_x,p_y)^T$ and write the equations of motion associated to
\equ{ham2rq} as
$$
\dot{\xi}=J\nabla H_0^{(q)}=M\xi\ ,
$$
where the matrices $J$ and $M$ are given by
$$
J=\begin{pmatrix}
0 & 0 & 1 & 0 \\
0 & 0&0&1 \\
-1& 0&0&0 \\
0 & -1 & 0 &0\\
\end{pmatrix}\ ,\qquad
M=\begin{pmatrix}
0 & 1 & 1 & 0 \\
-1 & 0&0&1 \\
2c_2 & 0&0&1 \\
0 & -c_2 & -1 &0\\
\end{pmatrix}\ .
$$
Next, after introducing the vector $\tilde \xi=(\tilde
x,\tilde y,\tilde z,\tilde p_x,\tilde p_y,\tilde p_z)$, one
performs a symplectic change of coordinates $\tilde \xi=C\,\xi$
with $C$ a real $4\times 4$ matrix (see \cite{CPS} for its
explicit expression), such that the transformed equations are
$$
{{d\tilde\xi}\over {dt}}=\widetilde M\ \tilde\xi
$$ 
with 
$$
\widetilde M=\begin{pmatrix}
\lambda_x & 0 & 0 & 0 \\
0 & 0&0&\omega_y \\
0 & 0&-\lambda_x &0 \\
0 & -\omega_y & 0 &0\\
\end{pmatrix}\ .
$$
In complete analogy with the circular case, we introduce the zero-order frequencies 
\beq{zfr}
\omega_y\equiv\sqrt{-\eta_1}, \quad \omega_z\equiv\sqrt{c_2}
\eeq  
and the rate
$$
\lambda_x\equiv\sqrt{\eta_2} 
$$
where
$$
\eta_1=\frac{c_2-2-\sqrt{9c_2^2-8c_2}}{2}, \qquad
\eta_2=\frac{c_2-2+\sqrt{9c_2^2-8c_2}}{2}\ . 
$$
As it is well known, 
since $c_2>1$, one obtains that $\eta_1<0$, $\eta_2>0$, showing
that the equilibrium point is of type saddle $\times$ center
$\times$ center and, as a matter of fact, the quadratic part $H_0^{(q)}$ can be written as 
$$
\tilde H_0^{(q)}(\tilde p_x,\tilde p_y,\tilde p_z,\tilde x,\tilde
y,\tilde z)=\lambda_x \tilde x\tilde p_x + \frac{\omega_y}{2}
(\tilde y^2+\tilde p_y^2)+\frac{\omega_z}{2}(\tilde z^2+\tilde
p_z^2)\ .
$$
Finally, we introduce complex coordinates through the canonical 
transformation 
\beqano
\tilde x&=&q_1\ , \qquad \qquad \ \ \tilde p_x=p_1\ ,\nonumber\\
\tilde y&=&\frac{q_2+ip_2}{\sqrt{2}}\ ,\qquad \tilde p_y=\frac{iq_2+p_2}{\sqrt{2}}\ ,\nonumber \\
\tilde z&=&\frac{q_3+ip_3}{\sqrt{2}}\ ,\qquad \tilde
p_z=\frac{iq_3+p_3}{\sqrt{2}}\ . 
\eeqano 
In the eccentric model,
such a transformation leads to write the quadratic part of the
Hamiltonian as 
\beqano
\H_{quadr}(p,q)&=&\lambda_x q_1p_1+i\omega_y q_2p_2+i\omega_z q_3p_3 - J_4\nonumber\\
&-& \sum_{n=1}^{\infty} \epsilon^n (-{\rm e})^n \cos^n f
\left[\frac12(x^2+y^2+z^2)+c_2(\mu) T_2 (x,y,z)\right]\ , 
\eeqano
where, in the square brackets, $x,y,z$ has to be intended as
expressed first in terms of $\tilde\xi$ and then in terms of $(p,q)$. We recognise that the zero-order
term in the extended phase space takes the diagonal form
\beq{H0D}
\H_0(p,q)=\lambda_x q_1p_1+i\omega_y q_2p_2+i\omega_z q_3p_3 + F
\eeq
and the starting Hamiltonian series can finally be written as in \eqref{ham3}
\beq{Hc}
\H (p,q,F,f)=\sum_{k\geq0} \epsilon^k \H_k(p,q,f)\ ,
\eeq
but expressed in the diagonalising coordinates.

\subsection{Center manifold reduction}\label{normf}
Perturbation theory is implemented here by constructing a 
resonant normal form generalising that obtained in the circular problem \cite{CPS,CCP}. 
It turns out that for typical values of the mass
parameter, the two elliptic frequencies $ \omega_y$ and $\omega_z $ in \eqref{zfr}
are very close to each other \cite{P19}; moreover, when extending the notion of frequencies to the elliptic problem, this fact is not substantially changed 
if the eccentricity is small. 
Therefore, in order to capture bifurcation 
phenomena related to the synchronous 1:1 resonance (e.g., \sl halo \rm orbits), we proceed to normalise with respect to the
\sl resonant \rm quadratic part 
\beqano
\H_{0}^{(res)}(p,q)&=&\lambda_x q_1p_1+i\omega_z (q_2p_2+
q_3p_3) - J_4\ .
\eeqano
Namely, we assume that the two elliptic frequencies $ \omega_y$ and $\omega_z $ are equal 
and we shift their small difference into
a first-order term
\beqano
\sum_{n=1}^{\infty} \epsilon^n \H_n &=& i\epsilon\delta q_2p_2 + \\
&-& \sum_{n=1}^{\infty} \epsilon^n (-{\rm e})^n \cos^n f \left[\frac12(x^2+y^2+z^2)+c_2(\mu) T_2 (x,y,z)\right] + \\
&-& \sum_{n=0}^{\infty} \epsilon^n (-{\rm e})^n \cos^n f
\sum_{j=1}^{\infty} \epsilon^j c_{j+2}(\mu) T_{j+2} (x,y,z)\ ,
\eeqano
by means of the \sl detuning \rm parameter $\delta$ defined as
$$
\delta \equiv \omega_y-\omega_z\ .
$$

In the end, we perform the \sl center manifold reduction \rm
(\cite{JM}), according to the following result.

\begin{proposition}\label{pro:cm}
Consider the Hamiltonian $\H (p,q,J_4,f)$ in \equ{Hc}. There exists a canonical
transformation $\C:\real^6\rightarrow\real^6$ with $\C(p,q)=(P,Q)$,
such that $\H$ is transformed to order $N\in\nat$ into the `normal form'
\beqano
K(P,Q,J_4)&=&\lambda_x Q_1P_1+i\omega_z ( Q_2P_2+  Q_3P_3)+i\epsilon\delta Q_2P_2 + F \nonumber\\
&+& \sum_{n=1}^N \epsilon^n
K_{n}(Q_1P_1,P_2,P_3,Q_2,Q_3)+R_{N+1}(P,Q)\ ,
\eeqano
where the polynomials $K_{n}$ depend on $Q_1,P_1$ only through their product $Q_1P_1$,
while $R_{N+1}(P,Q)$ is the remainder function of degree $N+1$, which
can depend on $Q_1$, $P_1$, separately.
\end{proposition}

We notice that the terms $K_{n}, n > 0$, satisfy the following properties:
\begin{itemize}
    \item[$(i)$] $K_{n}$ are polynomials of degree $n+2$ in $(P,Q)$;
\item[$(ii)$] the terms $K_{n}$ satisfy the following \sl normal form equation: \rm
$$
\{K_{n},H_0^{(res)}\}\big\vert_{(P,Q)}=0\qquad \forall
\, n > 0\ , 
$$
namely they are combinations of monomials in the \sl kernel \rm of the linear Hamiltonian operator associated to $H_0^{(res)}$;
    \item[$(iii)$] the coefficients of the terms $K_{n}$ depend on e, but all terms depending on $f$ are eliminated.
\end{itemize}

Summarising, the resonant normal form allows us to obtain three 
relevant results:

1. The normal form depends on the first
conjugate pair only through powers of the product $Q_1P_1$;
the hyperbolic dynamics is `killed' and the center manifold 
reduction is henceforth obtained by choosing initial conditions such that $Q_1=0,P_1=0$
(compare with Section~\ref{sec:NF1}). 

2. We pass from a non-autonomous to an autonomous integrable 
3-DOF hamiltonian system. The dependence on $f$, which is removed 
in the normalising transformation, is maintained in the generating
functions which are used to construct approximating solutions for
the orbits explicitly depending on the true anomaly of the primaries. 

3. The resonant terms allows us to describe bifurcation phenomena associated 
to the (almost) 1:1 commensurability of the linear frequencies $\omega_y$ and $\omega_z$.

\section{The resonant normal form}\label{sec:NF1}

By exploiting a recursion procedure based on the Lie-transform approach \cite{CPS}, according to Proposition 1, we construct the normal form, 
\beqa{NFS}
K(P,Q)
&=&\sum_{n\geq0}^N K_{2n}(P,Q) = \nonumber\\
&=& \lambda_x Q_1P_1+i\omega_z (Q_2P_2 + Q_3P_3)+i \delta Q_2P_2 \nonumber\\
&+& \sum_{n=1}^N K_{2n}(Q_1P_1,P_2,P_3,Q_2,Q_3)\ ,
\eeqa
which actually is of order $2N$ since, by symmetry, odd-degree terms
do not appear. We have ignored the higher-order remainder and put equal to one the book-keeping parameter. By construction, \eqref{NFS} does not depend anymore on $f$ so that $F$ is formally conserved and can be ignored from hereinafter. 

At each step of the recursion, we get the generating
functions $\chi_{n}$ as 
$$
\chi_{n}(P,Q,f), \qquad n=1,2,\dots,2N-1\ ,
$$
where each term is of degree $n+2$ in $(P,Q)$ and possibly trigonometric in $ \cos n f $. Following \cite{CCP}, the normal
form is expressed in the most compact way by using action-angle
variables defined through \beqa{QP}
Q_1&=&\sqrt{J_x}e^{\theta_x} \nonumber\\
Q_2&=&\sqrt{J_y}(\sin \theta_y -i\cos \theta_y)=-i\sqrt{J_y}e^{i\theta_y} \nonumber\\
Q_3&=&\sqrt{J_z}(\sin \theta_z -i\cos \theta_z)=-i\sqrt{J_z}e^{i\theta_z} \nonumber\\
P_1&=&\sqrt{J_x}e^{-\theta_x} \nonumber\\
P_2&=&\sqrt{J_y}(\cos \theta_y -i\sin \theta_y)=\sqrt{J_y}e^{-i\theta_y}  \nonumber\\
P_3&=&\sqrt{J_z}(\cos \theta_z -i\sin
\theta_z)=\sqrt{J_z}e^{-i\theta_z}\ .\nonumber 
\eeqa 
As said above $J_x=Q_1P_1$ is a formal conserved quantity of the normal form.  
The motion on
the \sl center manifold \rm $J_x=0$ is then described by
an effective 2-DOF Hamiltonian of the form
\beq{nfcm} 
K^{(CM)}(J_y,J_z,\theta_y,\theta_z)= \Omega_y J_y +
\Omega_z J_z + \sum_{n=1}^N K_{2n}(J_y,J_z,\theta_y-\theta_z)\ , 
\eeq 
where the \sl new \rm elliptic frequencies can be expanded
as powers of the eccentricity up to a given order $0\leq M \leq N$:
\beq{freqs}
\Omega_y = \sum_{n=0}^{M} \Omega_{y,n} \ {\rm e}^n\ ,\qquad \Omega_z =
\sum_{n=0}^{M} \Omega_{z,n} \ {\rm e}^n\ 
\eeq
with
$$
\Omega_{y,0} = \omega_{y}, \quad \Omega_{z,0} = \omega_{z} \ .
$$

For further reference, the terms up to third order in the actions are
\beqa{nfcm2} 
K_2 &=& \alpha J_y^2 + \beta J_z^2+J_y J_z(\sigma+2 \tau
\cos 2(\theta_y-\theta_z)), \label{nfcm22}\\ 
K_4 &=&  \alpha_1 J_y^3 + \beta_1 J_z^3 + \sigma_1 J_y^2 J_z + \sigma_2 J_y J_z^2 
+2 (\tau_1 J_y^2 J_z + \tau_2 J_y J_z^2) \cos 2(\theta_y-\theta_z)\label{nfcm24}
\eeqa
for suitable constant coefficients $\alpha$, $\beta$, $\sigma$, $\tau$, 
$\alpha_1$, $\beta_1$, $\sigma_1$, $\sigma_2$, $\tau_1$, $\tau_2$.
The coefficients of these expansions depend not only on the
mass-ratio $\mu$, but also on the eccentricity of the motion of
the primaries as terms of the series in ${\rm e}^{2n}$. They progressively appear 
as much as higher is the truncation order of the normal form. For example, 
truncating at $K_4$, in $K_2$ there appear 
terms proportional to ${\rm e}^{2}$.

The Hamiltonian $K^{(CM)}$ defines an integrable dynamics, since
it depends just on the combination angle $\theta_y-\theta_z$: the additional 
formal integral of motion is
$$
\E=J_y+J_z\ .
$$
In principle,
any initial condition on the center manifold determines a known
invariant object, periodic or quasi-periodic. In particular, the
nonlinear \sl normal-modes \rm $J_y = {\rm const}, J_z=0$ and $J_y
=0, J_z= {\rm const},$ respectively produce \sl planar \rm and \sl
vertical \rm Lyapunov orbits. Resonant solutions satisfying
$\theta_y-\theta_z=\pm \pi/2$ and suitable values of $J_y, J_z,$
produce the \sl halo \rm families. Generic initial conditions
produce invariant (Lissajous) tori around each family.

These solutions can be mapped back to the original coordinates by
inverting the normalising transformation. To do this, let us
collectively denote the final coordinates corresponding to a given
invariant object as $W=(P,Q)$. Then, the original coordinates $w=(p,q)$
approximating the `actual' object, are given by a series of the form
$w=\sum_k \epsilon^k w_k$ with terms given by the sequence 
\beqa{solss}
w_0 &=&  W, \nonumber\\
w_1 &=& \{W,\chi_1\}, \nonumber\\
w_2 &=& \{W,\chi_2\} + \frac12 \{\{W,\chi_1\},\chi_1\}, \nonumber\\
\vdots &=& \vdots \nonumber\\
w_k &=& \{W,\chi_k\} + \Gamma_k(W,\chi_1,...,\chi_{k-1}) 
\eeqa
for a suitable function $\Gamma_k$ depending on the generating functions determined at the previous steps. 
The above formulae allow us to deduce the explicit time-dependence
of the solutions. Finally, to plot the solutions in the synodic
system, we have to invert the diagonalising and scaling
transformations (see, e.g., \cite{CCP}). It is worth to observe that either a normal 
mode or another family of periodic orbits in \sl generic position \rm 
of the normal form (\cite{SV}), in 
general, \sl do not \rm give a periodic orbit in the original coordinates. The 
dependence on the true anomaly injected in the terms \eqref{solss}
by the generating functions $\chi_k$ appears at the frequency 
of the primaries, which is in general not-commensurable with the period of the 
orbit in the normalising variables. True periodic orbits can appear as isolated 
objects only when initial conditions produce such a commensurability. We will provide some examples in Section \ref{sec:EMC}.

\section{Halo orbits}\label{sec:halo} 

The theory of bifurcation of the halo family proceeds in the same
way as treated in the circular problem; we refer to
\cite{CCP,CPS} for more results in the context of the circular
problem. Here, we limit the analysis to bifurcations at first-order which, for sake 
of understanding the relevant mechanism, is perfectly fit.

Hamilton's equations associated to $K^{(CM)}$ in
\equ{nfcm}-\equ{nfcm2} are given by 
\beqa{EMJ}
\dot{J_y}&=& 4 \tau J_y J_z \sin (2\theta_y-2\theta_z)\label{EMJ1}\\
\dot{J_z}&=& - 4 \tau J_y J_z \sin (2\theta_y-2\theta_z)\label{EMJ2}\\
\dot{\theta_y}&=&\Omega_y + [2\alpha J_y +\sigma J_z + 2 \tau J_z \cos(2\theta_y-2\theta_z)]\\
\dot{\theta_z}&=&\Omega_z + [2\beta J_z +\sigma J_y + 2 \tau J_y
\cos(2\theta_y-\theta_z)]\ .
\eeqa

From \equ{EMJ1}--\equ{EMJ2} we have that $\dot J_y+\dot J_z=0$, so that, in
this approximation, $J_y+J_z$ is a constant of motion. According
to \cite{PM14,CPS}, we introduce the following transformation of
coordinates 
\beqa{ER}
\E&=&J_y+J_z\nonumber\\
\R&=&J_y\nonumber\\
\nu &=&\theta_z\nonumber\\
\psi &=&\theta_y-\theta_z\ . 
\eeqa 
The transformed Hamiltonian $K^{(CM)}$, obtained implementing
\equ{ER}, reads as
\beq{KTR}
K^{(tr)}(\E,\R,F,\nu,\psi)=\Omega_z \E+\delta
\R+a \R^2+b \E^2+c \E \R+d(\R^2-\E \R)\cos(2\psi)\ ,
\eeq
where
\beqano
a&=&\alpha+\beta-\sigma\ ,\nonumber\\
b&=&\beta\nonumber\\
c&=&\sigma-2\beta\nonumber\\
d&=&-2 \tau\ .
\eeqano
From Hamilton's equations associated to \equ{KTR}, we have:
\beqano
\dot{\E}&=&0\nonumber  \\
\dot{\R}&=&2d \R (\R-\E)\sin(2\psi)\nonumber\\
\dot{\nu}&=&\Omega_z +2b\E+c\R-d\R\cos(2\psi)\nonumber\\
\dot{\psi}&=&\delta+2a\R+c\E+d(2\R-\E)\cos(2\psi)\ . 
\eeqano
Hence, we recognize that the equilibrium positions
$\psi=0$, $\pi$, $\pm{\pi\over 2}$ give $\dot\R=0$; according to
\cite{PM14,CPS}, $\psi=0,\pi$ are denoted \sl inclined \rm orbits
and $\psi=\pm{\pi\over 2}$ \sl loop \rm orbits. The loops are usually 
dubbed \sl halo orbits \rm in the terminology associated to the collinear points. The equilibria, 
say $\R_{eq}$, are constrained by the condition
$$ 
0 \le \R_{eq} \le \E\ , 
$$
so that threshold values of $\E$ will appear determining the 
bifurcation of the corresponding family. 
We denote by $\E_N$ a truncation of the series up to an integer
order $N$, say 
$$
\E_N=\sum_{k=1}^N C_k\ \delta^k\ 
$$
for suitable real coefficients $C_k$.
Then, we look for a relation on the bifurcating normal mode
between $\E$ and $E$, that is the energy associated with the
Hamiltonian \equ{ham1}. These equilibrium positions stem from
bifurcations from the normal modes when entering the synchronous
resonance, since the normal modes lose stability. As said above, normal modes
are the periodic orbits along a single axis, which
correspond to the solutions of $\R=\E$ for motions along the
$y$-axis and $\R=0$ for motions along the $z$-axis. Such normal
modes represent, respectively, an approximation of the planar and
vertical Lyapunov periodic orbits.

A computation to first order in $\delta$ of the energy
level $E_1$ corresponding to a bifurcation to halo orbits is given
by
$$
E_1= \Omega_z \E_1 = \Omega_z C_1 {\delta}\ ,
$$
which, going back to the original coefficients, gives the
bifurcation value 
\beq{bif}
E_1=\frac{ \Omega_z \delta}{\sigma-2(\alpha+\tau)}\ . 
\eeq

\section{Example: $L_1$ in the Earth-Moon case}\label{sec:EMC}
In order to give a workable case, we provide figures 
for the Earth-Moon case ($\mu=0.012150586$)
with explicit expressions for the equilibrium point $L_1$. Moreover, 
we make use of our analytical theory, for the parameters of the 
Earth-Moon system, to find approximate initial conditions for
Lyapunov and halo type orbits. Several results are presented here, 
whereas more cumbersome formulas are provided in the Appendix.

\subsection{Normal-mode frequencies}
\begin{figure}[h]
\center
\includegraphics[height=5cm]{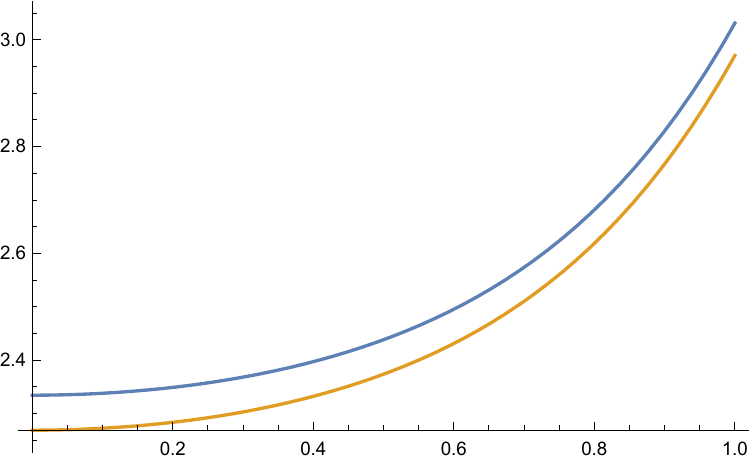}
\caption{$\Omega_y$ (blue) and $\Omega_z$ (yellow) for the point $L_1$ as a
function of the eccentricity in the Earth-Moon system, computed
using a normal form to order 2 and a series expansion of the 
eccentricity to order 8. 
}\label{fig:omegayz}
\end{figure}

We start with the horizontal and vertical frequencies \eqref{freqs} of the normal modes computed 
up to the power 8 in
the eccentricity:
\beqa{omyz}
\Omega_y &=&  2.33439 + 0.356732\, {\rm e}^2 + 0.200957\, {\rm e}^4 + 0.139319\, {\rm e}^6+ 0.106541\, {\rm e}^8\ ,\nonumber\\
\Omega_z &=&  2.26883 + 0.360261\, {\rm e}^2 + 0.201856\, {\rm
e}^4 + 0.139720\, {\rm e}^6+ 0.106766\, {\rm e}^8\ .
\eeqa
It is worth to observe that these expansions coincide term by term with those obtained in \cite{Holi} using the Poincar\'e-Lindstedt method. Figure~\ref{fig:omegayz} shows the graphs of $\Omega_y$ and $\Omega_z$ for the point $L_1$
as a function of the eccentricity, when computing a normal form to
the order $N=2$ and as a series expansion in the eccentricity to the order $M=8$.

\subsection{Planar Lyapunov orbits}\label{PLO}

By inserting the generating functions \equ{chi1}--\equ{chi2} 
given in the Appendix in the sequence of transformation \eqref{solss}, 
we get the explicit expressions \equ{xf}--\equ{zf} for the Cartesian coordinates. Examples of use of these series are 
provided in this and the following subsections.

\begin{figure}[h]
\center
\includegraphics[height=7.5cm]{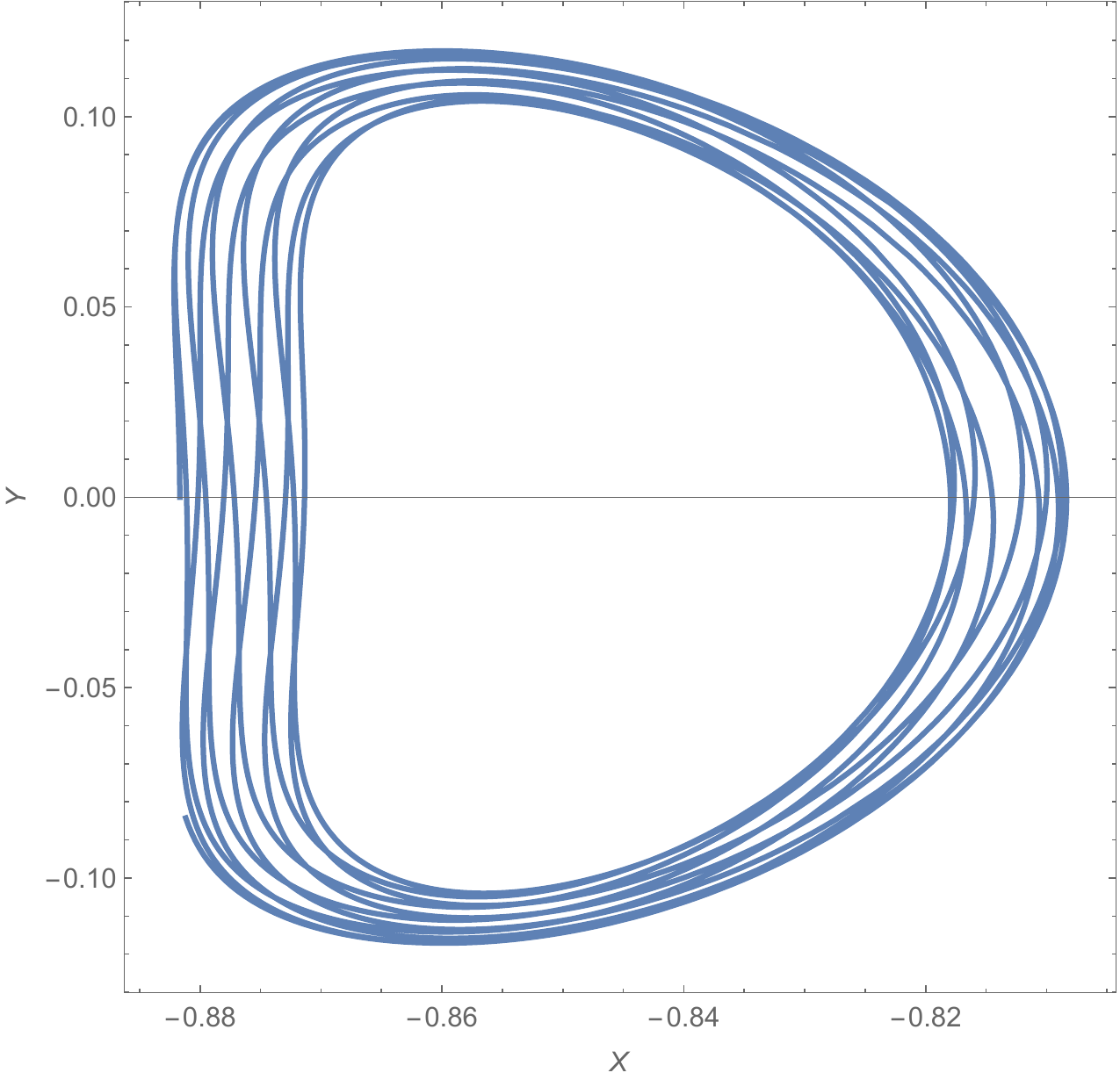}\\
\caption{A planar Lyapunov orbit around $L_1$ in the Earth-Moon
($\mu=0.012$) system with ${\rm e}=0.2$. Units are in terms of the Earth-Moon distance.
}\label{FR1}
\end{figure}

A normal form solution of the form $J_y = {\rm const}, J_z=0$
produces a planar Lyapunov orbit. Viceversa, a solution $J_y = 0, J_z= {\rm const}$ 
gives a \sl vertical \rm orbit. In the final normalising
variables, they correspond to nonlinear oscillations with
frequency 
\beqa{ffr1} 
\kappa_y&=&\Omega_y+2\alpha J_y+3\alpha_1 J_y^2\label{ffry}\\
\kappa_z&=&\Omega_z+2\beta J_z+3\beta_1 J_z^2\ . \label{ffrz}
\eeqa
As we mentioned before, for planar Lyapunov orbits, then $J_z=0$. 
Given $\kappa_y$, one determines $J_y$ through \equ{ffry}. 
When mapped back to the original space, since now the solution
explicitly depends on time, we see how a planar Lyapunov orbit
around $L_1$ of the \sl circular \rm problem is affected by the
elliptic motion of the primaries. As an example, using a big value
of the eccentricity in order to magnify its effect (${\rm e}=0.2$)
and taking $J_y=1$, we get a plot like that presented in Figure~\ref{FR1}.


\begin{figure}[h]
\center
\includegraphics[height=6.5cm]{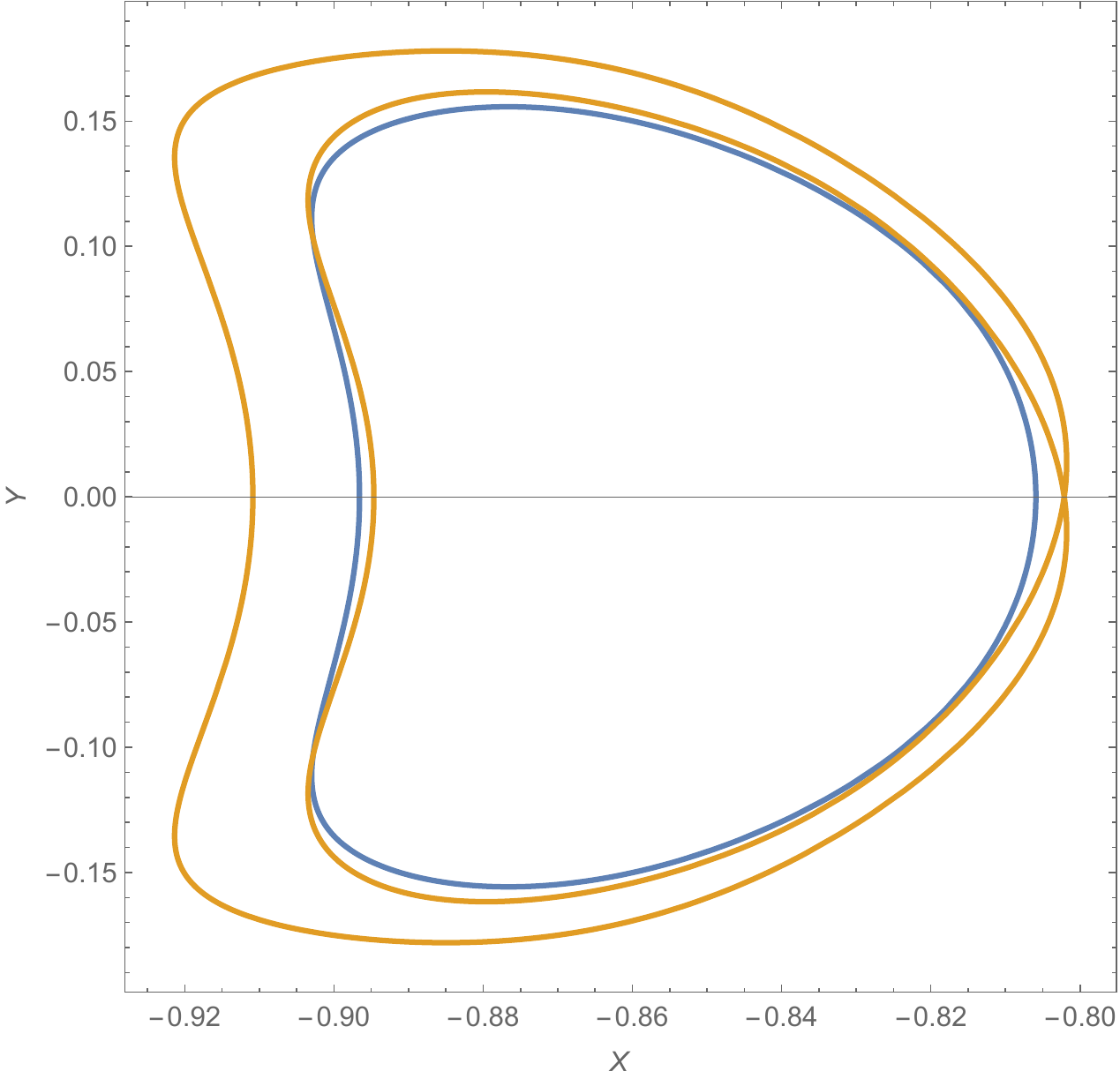}
\includegraphics[height=6.5cm]{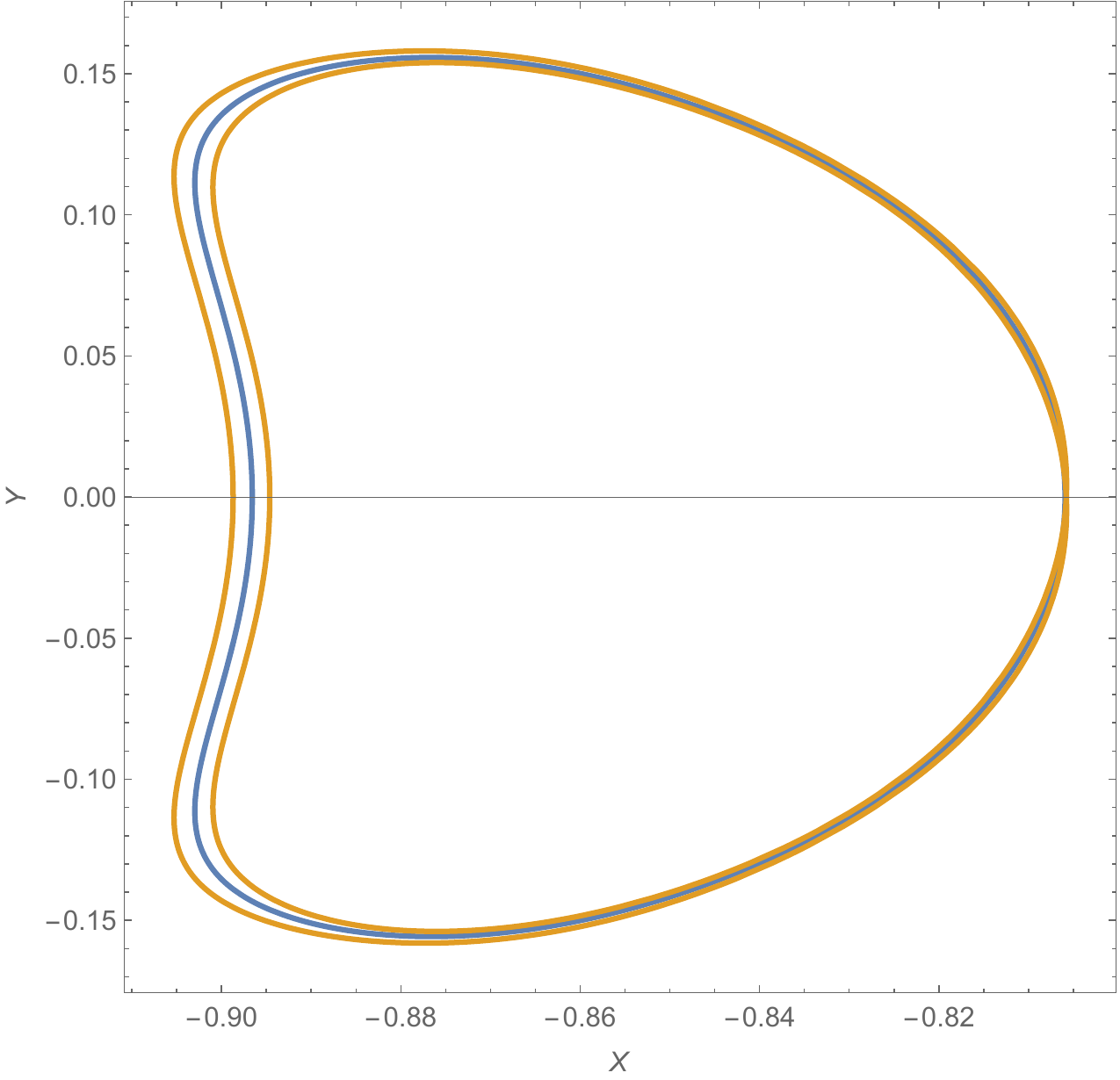}
\caption{Periodic planar Lyapunov orbits (in yellow) for $\kappa_y=$ 2:1 in
the Earth-Moon system: on the left we use ${\rm e}=0.2$ in order to magnify the effect. 
On the right the true value ${\rm e}=0.0549$ is used. 
The two orbits are respectively produced by the amplitudes $J_y=1.075$ and $J_y=0.9287$.
In blue the corresponding orbit in the circular case.}\label{FR2}
\end{figure}

In the eccentric case, it is possible to find \sl periodic \rm
Lyapunov orbits, provided that the frequency \eqref{ffr1}
satisfies a condition of the form $\kappa_y=m/n$ with $m,n$
suitable integers. This means that there is a resonance between
the orbital period and the synodic period. Once chosen a value for
the eccentricity, using expression \eqref{ffr1}, one obtains a
value of the amplitude which produces a certain commensurability.
In analogy with the example given in \cite{Holi}, where the
approximations have been found with the Lindstedt-Poincar\'e
method, we look for the first-order amplitude producing the
resonance condition $\kappa_y=$ 2:1. By using \eqref{ffr1} with
${\rm e}=0.2$, we find $J_y=1.075$. In the actual Earth-Moon case (${\rm e}=0.0549$) 
we find instead $J_y=0.9287$. 
By constructing the
corresponding solutions in the synodic system, we obtain a fairly
good approximation of periodic Lyapunov orbits with the 2:1 ratio, as shown by comparing  
the right panel of Figure~\ref{FR2} with Figure~\ref{FR8}.

%

\subsection{Vertical Lyapunov orbits}
The procedure to find periodic vertical Lyapunov is the same as that presented in Section~\ref{PLO}. For vertical Lyapunov orbits, then $J_y=0$ and equation \equ{ffrz} for a given $\kappa_z$ allows us to determine $J_z$. Then, Cartesian coordinates $X$, $Y$, $Z$ can be found in the Appendix and plotted as a function of the true anomaly $f$. We notice that, in analogy to the circular case, motion occurs also along $X$ and $Y$, although the normal form prescribes only an oscillation in the vertical direction. 

\begin{figure}[h]
\center
\includegraphics[height=7cm]{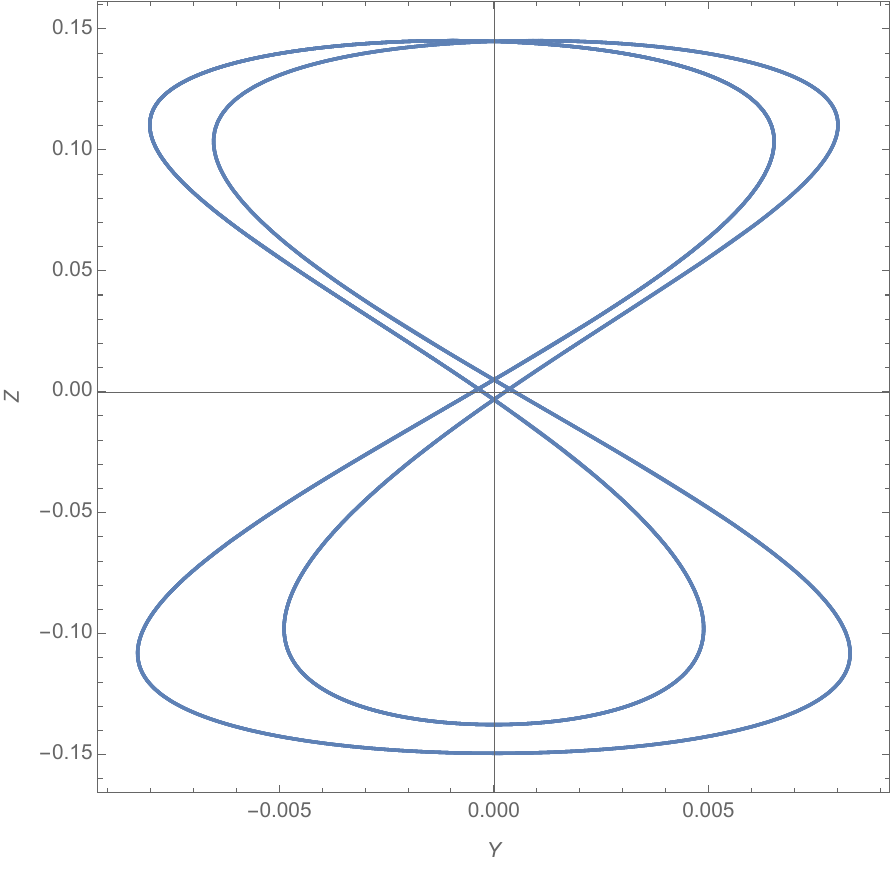}
\includegraphics[height=7cm]{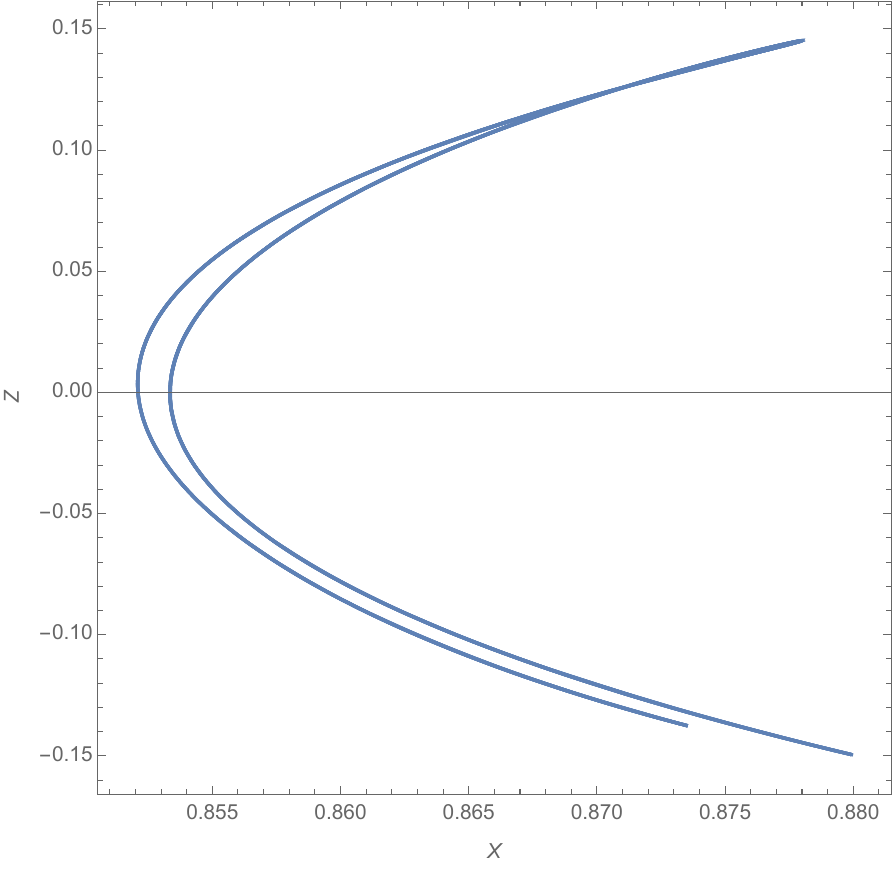}
\caption{Periodic vertical Lyapunov orbits for $\kappa_z=$ 2:1 in
the Earth-Moon ($\mu=0.012$) system but with ${\rm e}=0.2$. The
amplitude is $J_z=0.9775$. }\label{FR4}
\end{figure}

\begin{figure}[h]
\center
\includegraphics[height=7cm]{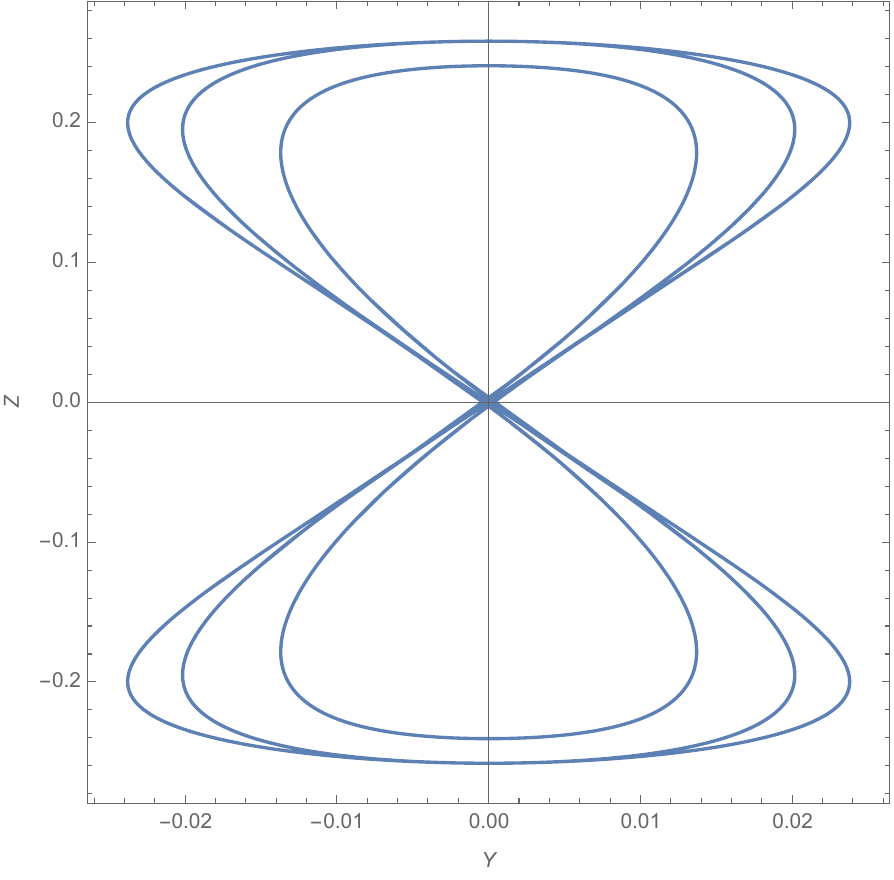}
\includegraphics[height=7cm]{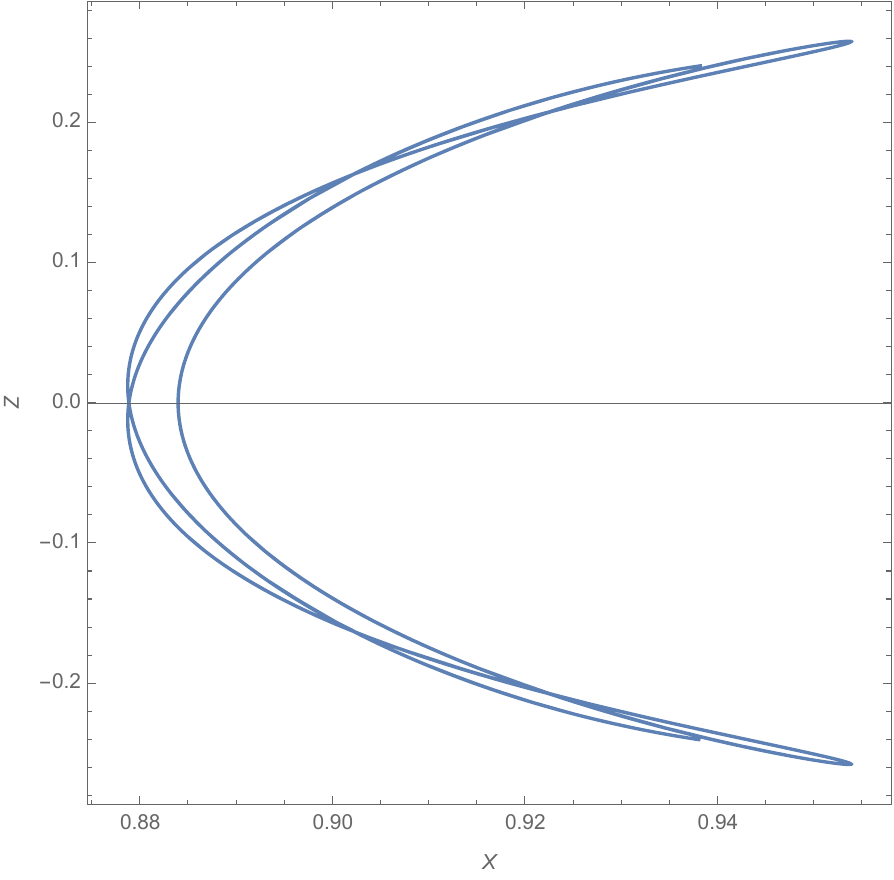}
\caption{Periodic vertical Lyapunov orbits for $\kappa_z=$ 3:2 in
the Earth-Moon ($\mu=0.012$) system but with ${\rm e}=0.2$. The
amplitude is $J_z=2.7030$. }\label{FR5}
\end{figure}

We report in Figure~\ref{FR4} an example with $\kappa_z=$ 2:1, while 
in Figure~\ref{FR5} we give an example with $\kappa_z=$ 3:2.


\subsection{Halo orbits}
In the case of the Earth-Moon system ($ {\rm e} = 0.0549$), using
(\ref{frL3num}-\ref{frL6num}), we get that the bifurcation value $E_1$ (see \equ{bif}) is given by 
$$
E_1= 0.30688 + 0.03221 {\rm e}^2 \simeq 0.30698\ ,
$$
which is slightly bigger than the value 0.306870 found on \cite{CPS} using a 2nd order resonant
perturbation theory.

In order to plot a halo orbit one can again exploit solutions \equ{xf}--\equ{zf}
of the Appendix. This procedure is clearly valid in the framework of a first-order
approximation, to keep the algebra as simple as possible.

One starts by choosing a value of $\E$ above the bifurcation value (4.9).
Then solves equation (4.7) for $\R_{eq}$ and finds:
$$ J_{y,H} = \frac{\E + \R_{eq}}2\ ,\qquad J_{z,H} = \frac{\E - \R_{eq}}2\ . $$

After that, one computes the frequency of the orbit:
$$ \Omega_H = \frac{\partial \K}{\partial \E} \bigg\vert_{\R=\R_{eq}} =
\omega_z
\frac{2\alpha (1 + 2 \beta \E) + 2\beta (1 + \delta) - (\sigma -
	2 \tau) (2 + \delta + (\sigma -
	2 \tau) \E)}{
	2 (\alpha + \beta - \sigma + 2 \tau)}
$$
and obtains
$$ \theta_y = \theta_z \pm \pi/2 = \Omega_H f . $$
Plus or minus in the angles give the \sl northern \rm and \sl southern \rm halo families.
Finally, one inserts $J_{y,H}, J_{z,H}, \theta_y, \theta_z$ into \equ{xf}--\equ{zf} and plots the
orbit.

\begin{figure}[h]
\center
\includegraphics[height=10cm]{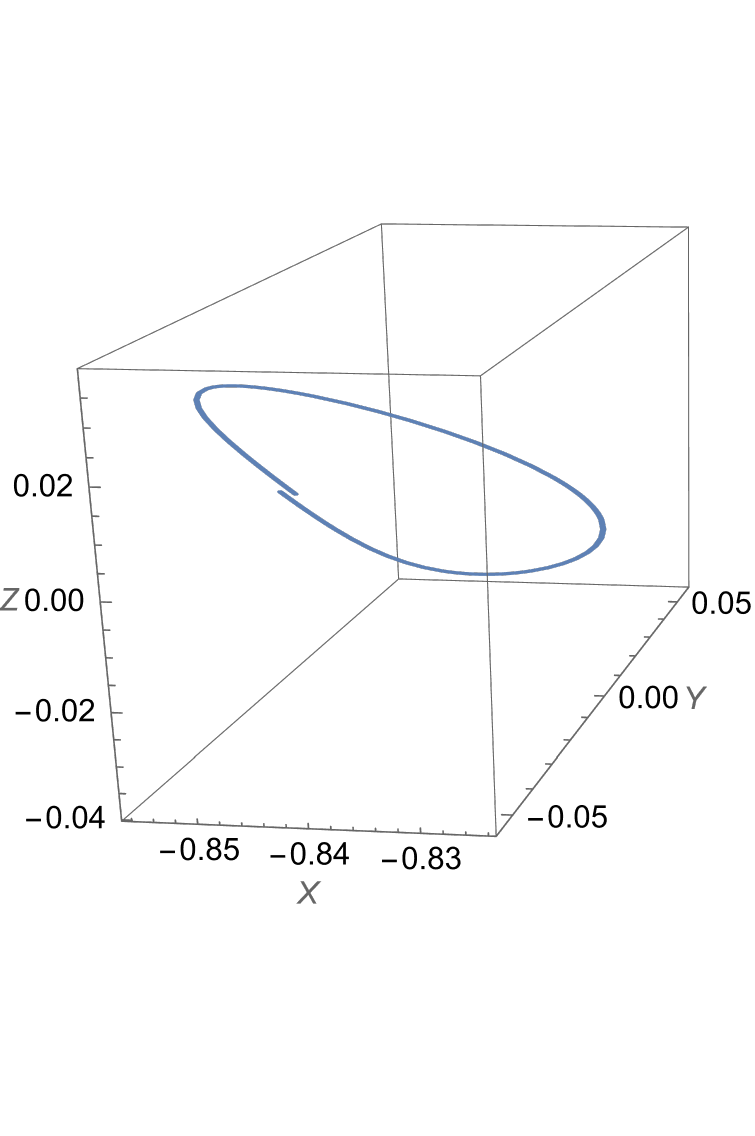}
\includegraphics[height=10cm]{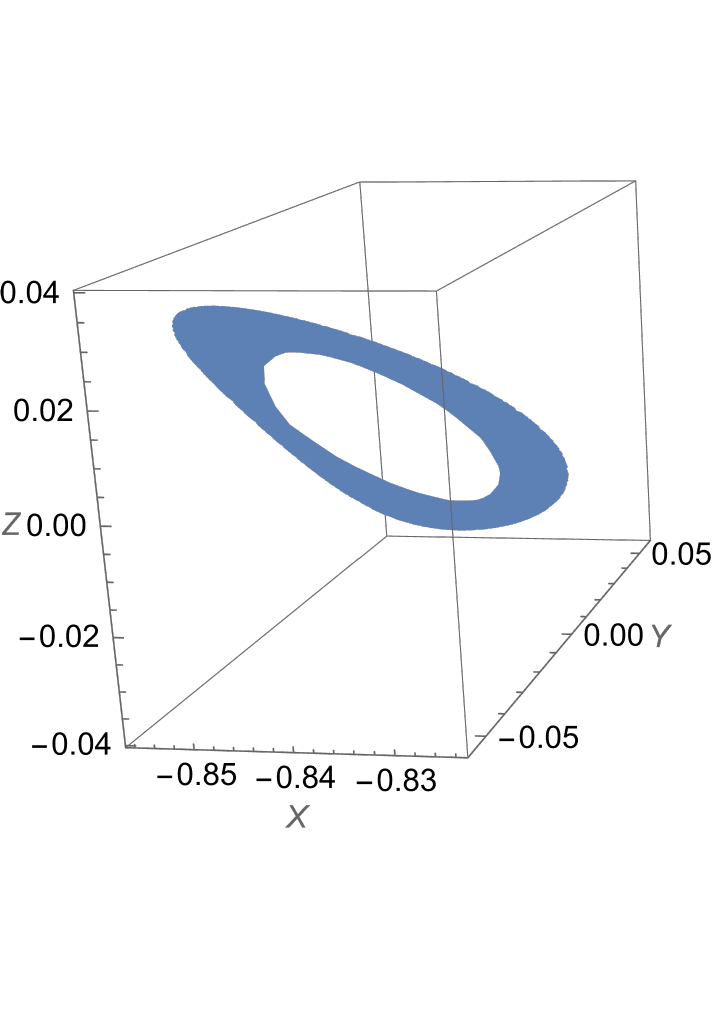}
\caption{Halo orbit in
the Earth-Moon system ($\mu=0.012$) around $L_1$: left, plotted for one-half synodic period; right, plotted for 100 synodic periods.}\label{FR6}
\end{figure}

In Figure~\ref{FR6} we show an example of a halo-orbit in the Earth-Moon system around $L_1$, 
plotted after one-half synodic period and after 100 synodic periods. 

\section{Validation by numerical simulations}\label{sec:numerical}

\begin{figure}[h]
\center
\includegraphics[height=6cm]{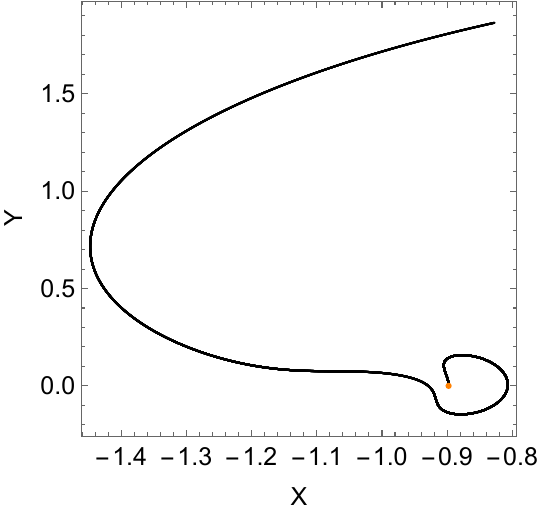}
\includegraphics[height=6cm]{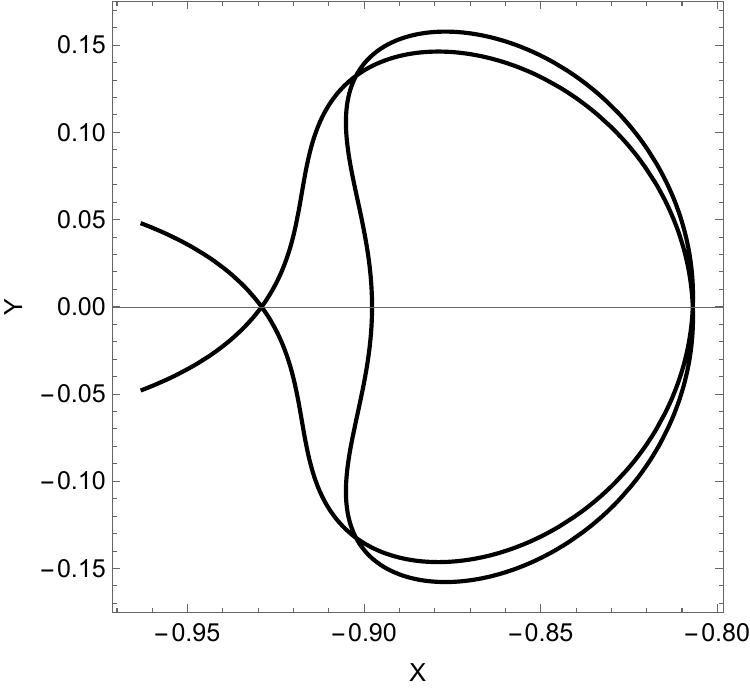}
	\caption{Same orbit as shown in Figure~\ref{FR2} for the Earth-Moon system with e=0.0549 based on direct numerical integration of \equ{eq1} with appropriate initial conditions computed using the analytic theory: one complete period $[0,T]$
	(left) and $[-T/2,T/2]$ (right) with $T=2\pi$.}\label{FR7}
\end{figure}

Special care needs to be taken to use the results from the previous sections
also in the original rotating-pulsating frame of reference. In this section we
highlight two major problems that may arise and also provide a numerical solution to
them.  We start by taking the approximate initial conditions (AICs), that were
obtained from the analytical theory, and use them as initial conditions for a
numerical integration of the system \equ{eq1}. As an example we provide the
orbit shown in Figure~\ref{FR7} using similar parameters and initial conditions
as for the 2:1 resonant planar orbit shown on the right of Figure~\ref{FR2}.  In
both panels of Figure~\ref{FR7}, we show the orbit in the plane $(X,Y)$ for one
synodic period $T$. From the analytical theory, we expect the orbit to perform two
closed loops during the integration time which we do not see on the left.
However, taking the time interval $[-T/2,T/2]$ the orbit approaches better the trajectory 
shown on the right of Figure~ \ref{FR2}. Still, we notice that 
$X'(-T/2)$ does not equal the value $X'(T/2)$ as we should expect for
a truly 2:1 resonant periodic planar orbit.
We think that the difference can
be explained by a combination of factors: i) the approximation error of the initial
conditions obtained from the analytical estimate and ii) errors in the numerical
integration.  Moreover, such errors are greatly amplified since the orbit is
unstable and the accumulation of the errors over the integration time grows
exponentially.  As a consequence, we conclude that a direct numerical integration of the AICs
using \equ{eq1} is not able to reproduce the correct orbit. 

To control the numerical problem ii), we use the solver NDSolve (Wolfram language, see
\cite{WL}) using arbitrary precision arithmetic (64 digits in our case). This 
choice is made to minimise the propagation of the errors due to numerical
instabilities when integrating unstable orbits in nonlinear dynamical systems.
Several test simulations, together with the analysis of the Jacobian of
\equ{eq1} along the solution vector, have been made to guarantee an accuracy of
the numerical integration up to 16-20 digits (using for safety a total of 64 digits).
Moreover, we checked that the Hamiltonian in the extended phase space is conserved
up to machine precision. As a result, we claim to have a reasonable control over the
errors stemming from problem ii). 

It remains to tackle the problem i). As it
turns out, even an eighth order normal form does not allow us to obtain a sufficient
number of digits in the AIC to overcome problem i).

However, we propose a solution that implements a numerical refinement of the
approximate analytical initial conditions.  
To this end, let $S=(X,Y,Z, X', Y', Z')$ be the state vector of
the orbit and let $S=S(t)$ be the state vector at a given time $t$. The orbit shown
in Figure~\ref{FR7} is based on the initial condition $S_0^{(0)}=S(t_0=0)$ where
$S_0^{(0)}=(X_0, Y_0, Z_0, X_0', Y_0', Z_0')$ denotes the AIC from
the normal form. Making use of the construction of the AIC, we can deduce the
following.  Since $S^{(0)}$ has been constructed to obtain a planar 2:1
resonant periodic orbit, we also have the property $S(T)=S(0)$ (called P1) with
given period $T$. Since the periodicity does not depend on the specific choice
of $f_0$ we also have $S(-T/2)=S(T/2)$ (called P2). Looking at the right
panel of Figure~\ref{FR7}, we see that property P2 is nearly fulfilled. Let us
denote by $\delta S_0^{(0)}$ the initial error of the approximate initial
condition labeled $S_0^{(0)}$ and define
$$
\delta S_0^{(0)} = S^{(0)}(-T/2) - S^{(0)}(T/2)\ , 
$$
where we used the upper-script $(0)$ in $S$ to indicate that the solution
$S=S(t)$ has been obtained using the initial condition $S_0^{(0)}$. The aim of the
numerical refinement method is to find better approximate initial conditions
$S_0^{(k+1)}$ such that $|\delta S_0^{(k+1)}|<|\delta S_0^{(k)}|$ with
$k=0,1,\dots$, where
\beqano
&&S_0^{(k)} = (X_0^{(k)}, Y_0^{(k)}, Z_0^{(k)}, X_0'^{(k)}, Y_0'^{(k)}, Z_0'^{(k)}) \\
&&\delta S_0^{(k)} = S^{(k)}(-T/2) - S^{(k)}(T/2)
\eeqano
and where the superscript $(k)$ denotes the $k$-th step of the numerical refinement. Thus,
we aim to minimize $|\delta S_0^{(k)}|$ up to $O(\varepsilon)$, say 
\beq{err2}
|\delta S_0^{(k)}| < \varepsilon ,
\eeq
for a given threshold $\varepsilon\in\mathbb{R}_+$.  In the case of a planar 2:1
resonant periodic orbit (right panel of Figure~\ref{FR2}), the approximate initial condition
has been constructed by making use of the symmetry of the periodic orbit that leads
to initial conditions of a specific form:
$$
S_0^{(0)} = (X_0^{(0)}, 0, 0, 0, Y_0'^{(0)}, 0) \ .
$$
In other words, from the geometry of the planar orbit, 4 out of 6 components vanish exactly:
\beq{err2b}
Y_0^{(k)} = Z_0^{(k)} = X_0'^{(k)} =  Z_0'^{(k)} = 0 \ ,\qquad k\in\integer_+\ .
\eeq
It thus remains to refine $(X_0^{(k)}, Y_0'^{(k)})$. A further reduction of
the complexity of the problem can be achieved as follows. In the circular problem the
so-called Jacobi constant $C_J$ is a conserved quantity along the orbit. Thus,
$C_J=C_J(X_0^{(k)}, Y_0'^{(k)})$ may be used to relate $X_0^{(k)}$ and
$Y_0'^{(k)}$ to define $\delta S_0^{(k)}$ to be a scalar function of a
single argument $x$ only: $\delta S_0^{(k)} = \delta S_0^{(k)} (x)$ (e.g., $x=X_0^{(k)}$). With this
we ensure that we converge to an orbit close to its initial Jacobi constant
(and also a related energy level).  In the elliptic case, however, $C_J=C_J(t)$
becomes a time dependent quantity.  Still, $C_J(0)$ relates $X_0^{(k)}$ with
$Y_0'^{(k)}$ at the initial time $t_0$. We make use of this relation, also in the 
elliptic problem, to define
the error function $F_k$ of a single variable $x$ to be a scalar function: 
\beq{err5}
F_k(x) = |\delta S_0^{(k)}(x)| \ .
\eeq
Thus, the problem to obtain refined initial conditions that fulfill \equ{err2} 
reduces to minimise \equ{err5} with respect to $x=X_0^{(k)}$. 

\begin{figure}[h]
\center
\includegraphics[height=7cm]{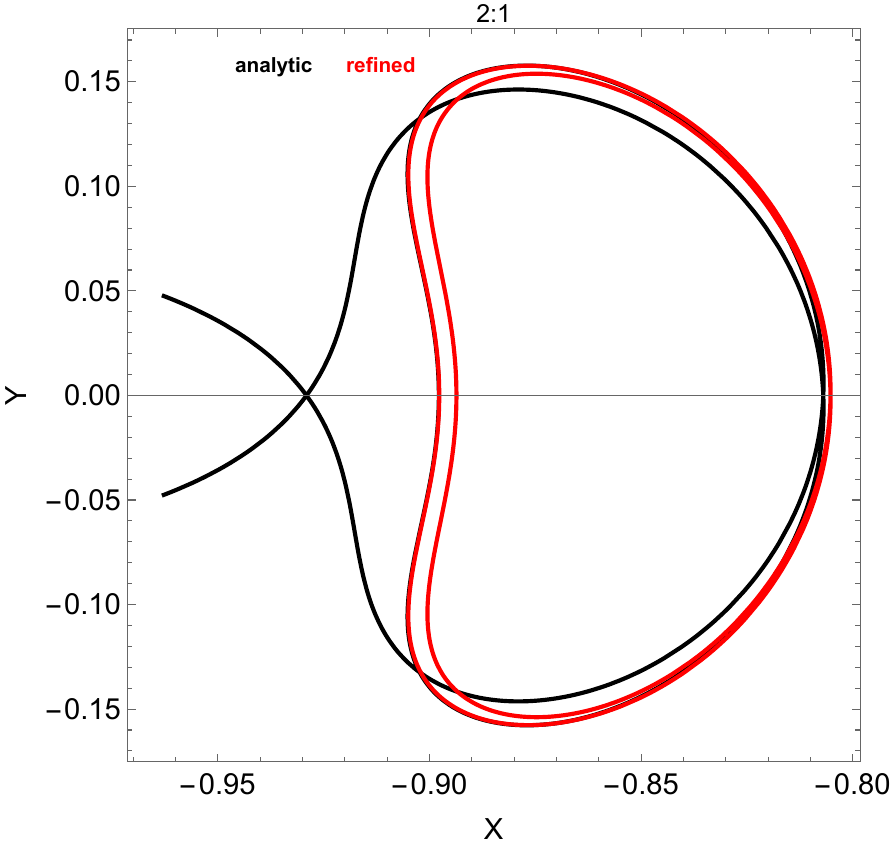}
	\caption{Same orbit as shown in Figure~\ref{FR7}. The black curve
	indicates the orbit obtain from approximate analytical initial
	conditions, the red curve gives the same orbit based on the refined
	approximate initial conditions (same integration time $T\simeq2\pi$).
	}\label{FR8}
\end{figure}

In practice, we integrate numerically the equations of motion \equ{eq1} using
$S_0^{(k)}$, varying $X_0^{(k)}$ until $F_k(X_0^{(k)})=O(\varepsilon)$, i.e. to
meet the requirement provided in \equ{err2}. We therefore notice the presence of an
additional parameter in the problem that we need to determine, which is the
integration time $T/2$.  From the condition for the orbit to be in $2:1$
resonance with the synodic period we conclude $T=2\pi$. However, as it turns
out \equ{err2b} is not necessarily fulfilled at $f = \pm T/2$ for approximate
initial conditions $S_0^{(k)}$. A solution to this problem can be found by
introducing the integration time itself, let us say $T_{1/2}^{(k)}$, as a free
parameter to the problem. At each step of refinement of the initial conditions
$S_0^{(k)}$ we determine $T_{1/2}^{(k)}$ by making use of the condition
\equ{err2b} itself.  In the planar case, we numerically integrate the equations
of motion using $S_0^{(k)}$ and determine $T_{1/2}^{(k)}$ from integration time
$f$ at which the condition $Y_0^{(k)}=0$ is fulfilled. We notice that for
approximate initial conditions, $X_0'^{(k)}$ does not necessarily vanish
when $Y_0^{(k)}=0$. Defining a second error function
\beq{err6}
G_k(x=X_0^{(k)}) = | X_0'^{(k)} | \quad {\rm when} \quad Y_0^{(k)} = 0 , 
\eeq
the correct resonant orbit is therefore determined by the requirement that both
conditions, \equ{err5} and \equ{err6}, vanish together. The implementation of
the procedure outlined above requires special care and comes with some solvable
technical difficulties: the numerical integrations based on \equ{eq1} need to
be very precise, and we need to check the conservation of the Hamiltonian in the
extended phase space. In addition, the numerical integrations need to be able to
determine very accurately the condition $Y_0^{(k)}=0$ to calculate
$G_k(X_0^{(k)})$. We notice that the convergence of the method does not only
depend on the initial choice $S_0^{(k)}$, but also on the choice of the solver
to minimise \equ{err5} and \equ{err6}.  In our implementation a combination of
a shooting method and Newton's method was used to ensure convergence during the
refinement process. 
We show in Figure~\ref{FR8} the orbit, mentioned at the beginning of the section, in the $(X, Y)$  plane
obtained by implementing the refinement of the analytical initial conditions.
The original orbit is shown in black (compare with Figure~\ref{FR7}, right
panel), while the numerical refined method is shown in red.  For the red orbit
$\varepsilon$ has been set equal to $10^{-32}$. This orbit should be compared
with the analytical approximation as shown in Figure~\ref{FR2}, right panel,
based on the normal form method, the agreement being very good. 

We remark that the method outlined above can be adapted to different kinds of 
orbit, i.e. resonant and non-resonant orbits, planar and spatial 
orbits in the circular and elliptic restricted three-body problem. 


\section*{Acknowledgements}
C.L. acknowledges the support from the Excellence Project 2023-2027 MatMod@TOV awarded
to the Department of Mathematics, University of Rome Tor Vergata, the project 
MIUR-PRIN 20178CJAB "New Frontiers of Celestial Mechanics: theory and applications", 
and the support of GNFM/INdAM.

G.P. acknowledges the support of INFN (Sezione di Roma2) and of GNFM/INdAM.


\section*{Appendix}

We provide explicit expressions for the normal form, the generating functions and
solutions for Lyapunov or halo orbit in terms of time (actually the true anomaly) for the case of the Earth-Moon system. We limit the results to order 4 ($N=2$). 
In order to show how the eccentricity of the primaries propagate into the results,  
we keep it in all formulae. The actual value of the Earth-Moon system is e=0.0549.

\subsection*{The normal form}
The terms linear in the actions are recovered by recalling \equ{nfcm}  
and looking at \equ{omyz}. 
For the coefficients appearing
in \eqref{nfcm22} and \eqref{nfcm24}, we get the following expressions: 
\beqa{calfas}
\alpha &=&  -0.162101380 -0.005834207  {\rm e}^2, \label{frL3num}\\
\beta &=&  -0.144882524 -0.006620655 {\rm e}^2, \label{frL4num}\\
\sigma &=&  -0.072614915 - 0.009081915  {\rm e}^2, \label{frL5num}\\
\tau &=& -0.23307061 - 0.002303888   {\rm e}^2 \label{frL6num} 
\eeqa
and
\beqano
\alpha_1 &=&  -0.01326986, \label{sfrL3num}\\
\beta_1 &=&  -0.008427193, \label{sfrL4num}\\
\sigma_1 &=&  -0.00294965, \label{sfrL5num}\\
\sigma_2 &=&  -0.0023065, \label{sfrL5num1}\\
\tau_1 &=& -0.03174666, \label{sfrL6num} \\
\tau_2 &=& -0.02757057. \label{sfrL6num1} 
\eeqano

\subsection*{The generating functions} 
The first two generating functions (when reduced to the center manifold) are 
\beqa{chi1}
\chi_1 &=& 0.12676787 J_y^{3/2} \sin \theta_y + 0.25911415 \sqrt{J_y} J_z \sin \theta_y - \nn
 && 0.052629037 J_y^{3/2} \sin 3 \theta_y + \nn
 && 0.13726658 \sqrt{J_y} J_z \sin(\theta_y - 2 \theta_z) - \nn 
 && 0.044009618 \sqrt{J_y} J_z \sin(\theta_y + 2 \theta_z) + \nn 
 && 0.11318452 {\rm e} J_y \sin (2 \theta_y - f) + 0.082529208 {\rm e} J_z \sin (2 \theta_z - f) + \nn 
 && 0.83480871 {\rm e} J_y \sin f + 0.91403775 {\rm e} J_z \sin f  + \nn 
 && 0.17488611 {\rm e} J_y \sin (2 \theta_y + f) + 0.12918669 {\rm e} J_z \sin (2 \theta_z + f)
 \eeqa
 and 
 \beqa{chi2}
\chi_2 &=& -0.086929879 {\rm e}^2 J_y \sin 2 \theta_y + 0.073423235 J_y^2 \sin 2 \theta_y + \nn
 && 0.082212129 J_y J_3 \sin 2 \theta_y - 0.026324747 J_y^2 \sin 4 \theta_y - \nn
 && 0.079393462 {\rm e}^2 J_3 \sin 2 \theta_z + 0.04150184 J_y J_3 \sin 2 \theta_z + \nn
 && 0.051258579 J_3^2 \sin 2 \theta_z - 0.0096648435 J_3^2 \sin 4 \theta_z - \nn
 && 0.033571224 J_y J_3 \sin (2 \theta_y + 2 \theta_z) - \nn
 && 0.075207468 {\rm e}^2 J_y \sin (2 \theta_y - 2f) - \nn 
 && 0.061131682 {\rm e}^2 J_3 \sin (2 \theta_z - 2f) - \nn
 && 0.012948487 {\rm e} J_y^{3/2} \sin (\theta_y - f) - \nn
 && 0.0096070536 {\rm e} \sqrt{J_y} J_3 \sin (\theta_y - f) + \nn
 && 0.038170957 {\rm e} J_y^{3/2} \sin (3\theta_y - f) + \nn
 && 0.11871541 {\rm e} \sqrt{J_y} J_3 \sin (\theta_y - 2 \theta_z- f) + \nn 
 && 0.033351938 {\rm e} \sqrt{J_y} J_3 \sin (\theta_y + 2 \theta_z- f) - \nn
 && 0.17836577 {\rm e}^2 J_y \sin 2f - 0.18013036 {\rm e}^2 J_3 \sin 2f + \nn 
 && 0.085110721 {\rm e} J_y^{3/2} \sin (\theta_y + f)+ \nn
 && 0.19449047 {\rm e} \sqrt{J_y} J_3 \sin (\theta_y + f) - \nn
 && 0.035544223 {\rm e} J_y^{3/2} \sin (3\theta_y + f) - \nn
 && 0.00060004321 {\rm e} \sqrt{J_y} J_3 \sin (\theta_y - 2 \theta_z + f) - \nn 
 && 0.02959452 {\rm e} \sqrt{J_y} J_3 \sin (\theta_y + 2 \theta_z + f) + \nn
 && 0.035853824 {\rm e}^2 J_y \sin (2\theta_y + 2 f) + \nn
 && 0.015525144 {\rm e}^2 J_3 \sin (2\theta_z + 2 f)\ .
\eeqa
Note that all terms explicitly depending on the true anomaly $f$ vanish when e=0.

\subsection*{Orbits}
Here we report the analytic solutions in Cartesian coordinates in the synodic frame for 
orbits in the Earth-Moon case with, as above, arbitrary eccentricity of the primaries. 
These expressions can be used to plot approximate orbits, when the frequencies are 
computed as in \equ{ffr1} or its equivalent for $\kappa_z$ and substituted into
$$ \theta_y = \kappa_y (J_y, J_z) f, \quad \theta_z = \kappa_z (J_y, J_z) f, $$
once given values of the actions $J_y, J_z$ are chosen. As discussed above in 
Section~\ref{PLO} this is particularly useful when looking for approximations of periodic orbits. 
In the analytical plots of Section \ref{sec:EMC} we made the further simplifying assumption of 
approximating $f$ with the mean anomaly $\ell = n t$. The Cartesian coordinates $X$, $Y$, $Z$ are given by the following expressions: 

\beqa{xf}
X(f) &=& -0.83691513 - 0.028066789 J_y - 0.028684328 J_z - \nn
 && 0.042530591 \sqrt{J_y} \cos \theta_y + 0.011616838 {\rm e}^2 \sqrt{J_y} \cos \theta_y - \nn
 && 0.0072538664 J_y^{3/2} \cos \theta_y - 0.0080462385 \sqrt{J_y} J_z \cos \theta_y + \nn
 && 0.011473421 J_y \cos 2 \theta_y + 0.0030761579 J_y^{3/2} \cos 3 \theta_y + \nn
 && 0.0012147932 \sqrt{J_y} J_z \cos (\theta_y - 2 \theta_z) + 0.012054911 J_z \cos 2 \theta_z + \nn 
 && 0.002869961 \sqrt{J_y} J_z \cos (\theta_y + 2 \theta_z) - 
       0.0004757866 {\rm e}^2 \sqrt{J_y} \cos (\theta_y - 2 f) + \nn 
 && 0.0038943343 {\rm e} \sqrt{J_y} \cos (\theta_y - f) - 
       0.0082945643 {\rm e} J2 \cos (2 \theta_y - f) - \nn
 && 0.009169466 {\rm e} J_z \cos (2 \theta_z - f) - 
       0.0083767837 {\rm e} J_y \cos f - 0.010112763 {\rm e} J_z \cos f - \nn
 && 0.037106626 {\rm e} \sqrt{J_y} \cos (\theta_y + f) + 
       0.013195973 {\rm e} J_y \cos (2 \theta_y + f) + \nn
 && 0.015384742 {\rm e} J_z \cos (2 \theta_z + f) - 
       0.0082310824 {\rm e}^2 \sqrt{J_y} \cos (\theta_y + 2f)\ ,
\eeqa

\beqa{yf}
Y(f) &=& 0.15253593 \sqrt{J_y} \sin \theta_y - 0.038529801 {\rm e}^2 \sqrt{J_y} \sin \theta_y + \nn
 && 0.011312977 J_y^{3/2} \sin \theta_y + 0.0030420108 \sqrt{J_y} J_z \sin \theta_y + \nn
 && 0.011580803 J_y \cos \theta_y \sin \theta_y + 0.0040746736 J_y^{3/2} \sin 3\theta_y - \nn
 && 0.0079110266 \sqrt{J_y} J_z \sin (\theta_y - 2 \theta_z) - \nn
 && 0.01330702 J_z \cos \theta_z \sin \theta_z + \nn
 && 0.001691008 \sqrt{J_y} J_z \sin (\theta_y + 2 \theta_z) + \nn
 && 0.018042238 {\rm e}^2 \sqrt{J_y} \sin (\theta_y - 2f) - \nn
 && 0.049106567 {\rm e} \sqrt{J_y} \sin (\theta_y - f) - \nn
 && 0.0039015441 {\rm e} J_y \sin (2 \theta_y - f) + \nn
 && 0.004349965 {\rm e} J_z \sin (2 \theta_z - f) + 0.0079996739 {\rm e} J_y \sin f + \nn 
 && 0.0064257087 {\rm e} J_z \sin f + 0.091803411 {\rm e} \sqrt{J_y} \sin (\theta_y + f) + \nn 
 && 0.010214749 {\rm e} J_y \sin (2 \theta_y + f) - \nn
 && 0.011359961 {\rm e} J_z \sin (2 \theta_z + f) + \nn
 && 0.006618954 {\rm e}^2 \sqrt{J_y} \sin (\theta_y + 2f)\ ,
\eeqa

\beqa{zf}
Z(f) &=& -0.037811696 \sqrt{J_y} \sqrt{J_z} \sin (\theta_y - \theta_z) - \nn
 && 0.005705595 J_y \sqrt{J_z} \sin (2 \theta_y - \theta_z) - \nn
 && 0.14171043 \sqrt{J_z} \sin \theta_z + 0.039184292 {\rm e}^2 \sqrt{J_z} \sin \theta_z - \nn 
 && 0.0013614348 J_y \sqrt{J_z} \sin \theta_z - 0.0068725065 J_z^{3/2} \sin \theta_z - \nn
 && 0.0027392183 J_z^{3/2} \sin (3 \theta_z) - \nn
 && 0.012122967 \sqrt{J_y} \sqrt{J_z} \sin (\theta_y + \theta_z) - \nn 
 && 0.0051006101 J_y \sqrt{J_z} \sin (2 \theta_y + \theta_z) - \nn 
 && 0.019654877 {\rm e}^2 \sqrt{J_z} \sin (\theta_z - 2f) - \nn
 && 0.019284097 {\rm e} \sqrt{J_y} \sqrt{J_z} \sin (\theta_y - \theta_z -f) + \nn 
 && 0.053069093 {\rm e} \sqrt{J_z} \sin (\theta_z - f) + \nn
 && 0.0081153711 {\rm e} \sqrt{J_y} \sqrt{J_z} \sin (\theta_y + \theta_z -f) - \nn 
 && 0.018253022 {\rm e} \sqrt{J_y} \sqrt{J_z} \sin (\theta_y - \theta_z +f) - \nn
 && 0.083071445 {\rm e} \sqrt{J_z} \sin (\theta_z + f) - \nn
 && 0.020585513 {\rm e} \sqrt{J_y} \sqrt{J_z} \sin (\theta_y + \theta_z +f) - \nn 
 && 0.0049915983 {\rm e}^2 \sqrt{J_z} \sin (\theta_z + 2f)\ .
\eeqa



\begin{thebibliography}{99}

\bibitem{Broucke}
Broucke, R., \it Stability of Periodic Orbits in the Elliptic,
Restricted Three-Body Problem, \rm AIAA Journal {\bf 7}, n. 6, 1003--1009 (1969)

\bibitem{AMPA}
Bucciarelli, S. Ceccaroni, M. Celletti, A. \& Pucacco, G. {\it
Qualitative and analytical results of the bifurcation thresholds
to halo orbits}, Annali di Matematica Pura e Applicata, {\bf 195},
489--512 (2016)

\bibitem{CCP}
Ceccaroni, M. Celletti, A.  \& Pucacco, G. \emph{Halo orbits
around the collinear points of the restricted three-body problem},
Physica D, {\bf 317}, 28--42 (2016)

\bibitem{Alebook}
Celletti, A. {\sl Stability and Chaos in Celestial Mechanics},
Springer-Verlag, Berlin; published in association with Praxis
Publishing Ltd., Chichester (2010)


\bibitem{CPS}
Celletti, A. Pucacco, G. Stella, D. {\sl Lissajous and Halo orbits
in the restricted three-body problem}, J. Nonlinear Science {\bf
25}, Issue 2, 343--370 (2015)

\bibitem{Conley}
Conley, C.~C., \sl Low energy transit orbits in the restricted three-body problem, \rm 
SIAM J. Appl. Math \textbf{16}, 732-746 (1968) 

\bibitem{FeLa}
Ferrari, F. Lavagna, M. \sl Periodic Motion Around Libration 
Points in the Elliptic Restricted Three-Body Problem
\rm
Nonlinear Dynamics, {\bf 93}, 453--462 (2018)

\bibitem{GM}
G\'omez, G. Mondelo, J.M., \sl The dynamics around the collinear equilibrium points of the RTBP, \rm Phys. D
{\bf 157}, 283–321 (2001)

\bibitem{JM}
Jorba, \`A. Masdemont, J. \sl Dynamics in the center manifold of
the collinear points of the restricted three body problem, \rm
Physica D {\bf 132}, 189--213 (1999)

\bibitem{Holi}
Hou, X. Y. Liu, L. {\sl On motions around the collinear libration
points in the elliptic restricted three-body problem}, MNRAS {\bf
415}, 3552--3560 (2011)

\bibitem{kolomaro}
Koon, W.S. Lo. W.L. Marsden, J.E. Ross, S.D., 
\sl Dynamical Systems, the Three-Body Problem
and Space Mission Design, \rm 
Mission Design Book Online (2011) 

\bibitem{LXD}
Luo, T. Xu, M. Dong, Y. {\sl Natural formation flying on
quasi-halo orbits in the photogravitational circular restricted
three-body problem}, Acta Astronautica {\bf 149}, 35--46 (2018)

\bibitem{PeXu}
Peng, H. Shijie, X. {\sl Stability of two groups of multi-revolution 
elliptic halo orbits in the elliptic restricted three-body problem}, Celestial Mechanics and Dynamical Astronomy, {\bf
123}, 279--303 (2015)

\bibitem{P19}
Pucacco, G. {\it Structure of the centre manifold of the $L_1, L_2$ collinear libration points in the restricted three-body problem}, Celestial Mechanics and Dynamical Astronomy, {\bf 131}, 44 (2019)

\bibitem{PM14}
Pucacco, G. Marchesiello, A. {\sl An energy-momentum map for the
time-reversal symmetric 1:1 resonance with $\Z_2\times\Z_2$
symmetry}, Physica D \textbf{271}, 10--18 (2014)

\bibitem{Richardson}
Richardson, D.L., \sl Analytic construction of periodic orbits about the collinear points, \rm Celest. Mech. \textbf{22},
241–253 (1980)

\bibitem{SV}
Sanders, J. A. Verhulst, F. and Murdock, J. \emph{Averaging Methods in Nonlinear Dynamical Systems}, Springer-Verlag, Berlin Heidelberg (2007)

\bibitem{Szebehely}
Szebehely, V. \sl Theory of Orbits, \rm Academic Press, New York
and London (1967)

\bibitem{WL}
	Wolfram Mathematica (Version 13.2), Wolfram Research (2024)
\url{https://reference.wolfram.com/language/ref/NDSolve.html}

\end{thebibliography}
\end{document}